\documentclass[11pt]{article}
\usepackage[utf8]{inputenc}
\usepackage{amssymb,graphicx,mathtools,amsfonts}
\usepackage{amsmath}
\usepackage{relsize}
\usepackage{amsthm}
\usepackage{hyperref}
\usepackage{lmodern}
\usepackage[dvipsnames]{xcolor}
\usepackage{comment}
\usepackage{tikz-cd}
\usepackage{enumitem}

\usepackage{upgreek} 
\newcommand{\Jq}{J_q}                               
\newcommand{\Lpq}{\mathcal{L}_{p,q}}                          
\newcommand{\RN}{\mathbb{R}^N}

\DeclareMathOperator{\dist}{dist}

\usepackage{a4wide}
\usepackage{esint}
\numberwithin{equation}{section}
\graphicspath{{images/}}

\newtheorem{theorem}{Theorem}[section]
\newtheorem{lemma}[theorem]{Lemma}
\newtheorem{definition}[theorem]{Definition}
\newtheorem{proposition}[theorem]{Proposition}
\newtheorem{corollary}[theorem]{Corollary}
\newtheorem{remark}[theorem]{Remark}

\newcommand{\R}{\mathbb{R}}
\allowdisplaybreaks

\newcommand{\om} {\Omega}

\newcommand{\rnn}{\mathbb{R}^{N}}

\newcommand{\X}{\mathcal{X}_0} 

\usepackage{tikz}

 \DeclareMathOperator{\dive}{div}

\newcommand{\dmu}{\,d\mu_{s,q}}

\newcommand{\lam}{\lambda}
\newcommand{\ps}{p^*}
\usetikzlibrary{positioning}

\usepackage[normalem]{ulem}
\newcommand\magsout{\bgroup\markoverwith{\textcolor{magenta}{\rule[0.5ex]{4pt}{1.5pt}}}\ULon}

\providecommand{\subjclass}[2][2020]{%
  \par\medskip\noindent\textbf{#1 MSC.}\ #2\par}
\providecommand{\keywords}[1]{%
  \par\smallskip\noindent\textbf{Keywords.}\ #1\par}

\title{Pohozaev identity, nonexistence and existence for mixed local-nonlocal $(p,q)$-Laplacian}

\author{Shammi Malhotra\footnote{Department of Mathematics, Indian Institute of Technology Delhi, Hauz Khas New Delhi 110016,  India, shammi22malhotra@gmail.com}, Diksha Gupta\footnote{Department of Mathematical Sciences, Korea Advanced Institute of Science and Technology (KAIST), Daejeon 34141, Republic of Korea, dikshagupta1232@gmail.com}, and K. Sreenadh$^{*}$}

\date{}

\begin{document}
\maketitle
\begin{abstract}
We study the mixed local-nonlocal operator $\Lpq=-\Delta_p+(-\Delta)^s_q$ in $\RN$. First, we prove a Pohozaev identity for bounded, continuous weak solutions of $\Lpq u=f(u)$, where $f$ is any continuous function. We impose no growth condition on $f$. Combined with the Nehari identity, this gives nonexistence of nontrivial solutions for nonlinearities that are supercritical, or subcritical, with respect to both operators. For $f(t)=|t|^{r-2}t$, this rules out all $r$ outside the interval between $p^*$ and $q^*_s$. Second, we show that this interval is exactly where solutions exist. When $p^*\neq q^*_s$ and $r \in \left(\min\{p^*,q^*_s\}, \;\max\{p^*,q^*_s\}\right)$, we obtain a nonnegative radial solution by constrained minimization and a ground state by using the Pohozaev identity as a constraint. Finally, we show that the two solutions have the same energy.
\end{abstract}

\section{Introduction}
In this paper, we study Pohozaev identity, nonexistence and existence results for a class of problems involving the mixed local-nonlocal operator
\[
\Lpq u:=-\Delta_p u+(-\Delta)^s_q u ,
\]
where $-\Delta_p u=-\dive(|\nabla u|^{p-2}\nabla u)$ is the $p$-Laplacian and $(-\Delta)^s_q$ is the fractional $q$-Laplacian of order $s\in(0,1)$ defined as 
\[
(-\Delta)^s_q u(x)
:= 2\,\mathrm{P.V.}\int_{\mathbb{R}^N}
\frac{|u(x)-u(y)|^{q-2}\big(u(x)-u(y)\big)}{|x-y|^{N+sq}}\,dy.
\] 
Operators which combine a local and a nonlocal diffusion have been studied intensively in recent years, see for instance \cite{BDVV2022,GarainKinnunen2022} and the references therein. When $p\neq q$ the two parts of $\Lpq$ have different homogeneities, and in general they also scale differently, so several tools which are standard for each of the two operators separately are not directly available for $\Lpq$. 
We denote by
\[
p^*:=\frac{Np}{N-p},\qquad q^*_s:=\frac{Nq}{N-sq}
\]
the critical Sobolev exponents of the two operators. We are particularly interested in studying the following problem:
\begin{equation}\label{eq:Pr}
\Lpq u=|u|^{r-2}u\qquad\text{in }\RN,\quad  r>1. \tag{$\mathcal{P}_r$}
\end{equation}

In~\cite{ourPreviousPaper_isocriticalregimemixedlocalnonlocal}, the authors studied problem~\eqref{eq:Pr} in the isocritical case
\begin{equation}\label{eq:isocritical-intro}
  p^*=q^*_s,
\end{equation}
in which the two parts of $\Lpq$ have the same scaling, and proved that the problem admits a solution for $r=p^*=q^*_s$, along with other results concerning qualitative properties such as decay estimates and comparison principles. This raises a natural question that what happens for other exponents $r$? For the purely local and the purely nonlocal operators it is known that, in the whole of $\RN$, there are no nontrivial solutions of finite energy for exponents other than the critical one, see for instance \cite{berestycki_lions_1983, Liu_plapPohozaev, RosOtonSerra2014}, and we expect the same here. The standard tool for such nonexistence results is the Pohozaev identity, and this leads us to study the same for $\Lpq$. Its proof does not use the isocritical condition, so we work with general $p$ and $q$, under the following minimal assumptions
\begin{equation}\label{assump: standing}
N\ge2,\qquad 1<p<N,\qquad 0<s<1,\qquad q>1,\qquad sq<N. \tag{H}
\end{equation}
Thus our first result is a Pohozaev identity for the problem
\begin{equation}\label{eq:P}
 \Lpq U=f(u)\qquad\text{in }\RN, \tag{$\mathcal P$}
\end{equation}
where the nonlinearity $f$ satisfies
\begin{equation}\label{assump:f}
  f\in C(\R),\qquad F(t):=\int_0^t f(\tau)\,d\tau\quad\text{for }t\in\R. \tag{F}
\end{equation}
Thus $F\in C^1(\R)$, $F'=f$ and $F(0)=0$. Here and below, $\X(\RN)=D^{1,p}_0(\RN)\cap D^{s,q}_0(\RN)$ is the energy space, $[\,\cdot\,]_{1,p}$ and $[\,\cdot\,]_{s,q}$ are the seminorms of the two operators, and $\langle\cdot,\cdot\rangle_p$, $\langle\cdot,\cdot\rangle_{s,q}$ are the associated forms (refer to Section~\ref{sec:prelim}).
\begin{definition}\label{def:weak}
A function $u\in\X(\RN)$ is a weak solution of \eqref{eq:P} if
\begin{equation}\label{eq:weak}
\langle u,\varphi\rangle_{p}+\langle u,\varphi\rangle_{s,q}=\int_{\RN}f(u)\,\varphi\,dx
\tag{W}
\end{equation}
for every $\varphi\in\X(\RN)$ with $f(u)\varphi\in L^1(\RN)$.
\end{definition}
Since $f$ is continuous, $f(u)$ is measurable, and $f(u)\in L^\infty(\RN)$ whenever $u\in L^\infty(\RN)$. If $p\le sq$, then $p<q$, and we set
\begin{equation}\label{eq:lstar}
    l_*:=\frac{sq-p}{q-p}\in[0,s).
\end{equation}

\begin{theorem}[Pohozaev identity]\label{thm:main Pohozaev}
Assume \eqref{assump: standing} and \eqref{assump:f}. Let $u\in \X(\RN)\cap L^\infty(\RN)\cap C(\RN)$ be a weak solution of \eqref{eq:P} such that
\begin{equation}\label{eq:integrability}
  f(u)\in L^{p'}(\RN)\qquad\text{and}\qquad F(u)\in L^1(\RN). \tag{I}
\end{equation}
If $p\le sq$, assume in addition that $u\in C^{0,l}_{\mathrm{loc}}(\RN)$ for some $l\in(l_*,1]$, with $l_*$ as in \eqref{eq:lstar}.
Then
\begin{equation}\label{eq:poho}
  \frac{N-p}{p}\,[u]_{1,p}^p
  +\frac{N-sq}{q}\,[u]_{s,q}^q
  =N\int_{\RN}F(u)\,dx.
\end{equation}
\end{theorem}

We stress that no growth condition on $f$ is needed in Theorem~\ref{thm:main Pohozaev}. This is because when $u$ is bounded, only the behavior of $f$ near $0$ matters in \eqref{eq:integrability}, and the growth of $f$ at infinity plays no role. In fact, set $m_-:=\min\{\ps,q^*_s\}$, and assume that there are constants $a\ge\max\{p,m_-\}$, $\delta>0$ and $C>0$ such that
\begin{equation}\label{eq:f-near-zero}
  |f(t)|\le C|t|^{a-1}\qquad\text{for }|t|\le\delta .
\end{equation}
Then every $u\in\X(\RN)\cap L^\infty(\RN)$ satisfies \eqref{eq:integrability}, and also $f(u)u\in L^1(\RN)$. For instance, one admissible function is $f(t)=|t|^{r-2}t\,e^{t^2}$ with $r\ge\max\{p,m_-\}$, which grows faster than any power at infinity.
\medskip 

For the Laplacian, the identity goes back to Pohozaev \cite{pohozaev1965}, and it was used by Br\'ezis and Nirenberg \cite{brezis1983} to rule out solutions of critical problems on star-shaped domains, for solutions of class $C^2(\om)\cap C^1(\bar\om)$ with $\om$ bounded. Guedda and V\'eron \cite{gueddaVeron_plap} proved it for the $p$-Laplace equation on a bounded domain, for solutions in $W^{1,p}_0(\om)\cap L^\infty(\om)$, and excluded positive solutions of $-\Delta_pu=u^{\ps-1}$ on star-shaped domains. For the full space, this was done by Liu and Liu \cite{Liu_plapPohozaev}, for solutions in $W^{1,p}(\RN)$. For extremals of a general variational problem, Pucci and Serrin \cite{PucciSerrin_Pohozaev} established one identity which contains all of these: for a Lagrangian $L(x,u,\nabla u)$ and a vector field $h\in C^1$, every critical point $u\in C^2(\om)\cap C^1(\bar\om)$ satisfies a divergence identity, the case $h(x)=x$ being the classical one. Their proof differentiates the Lagrangian, which excludes the $p$-Laplacian for $1<p<2$. Later, Degiovanni, Musesti and Squassina \cite{DegiovaniSquassina_Pohozaev} extended their result to $C^1$ solutions, assuming instead that $L$ is strictly convex in the gradient.
\medskip

For nonlocal operators, modeled on the fractional Laplacian, proving such an identity on a bounded domain was a longstanding problem, since there is no pointwise chain rule to integrate by parts. It was resolved much later by Ros-Oton and Serra \cite{RosOtonSerra2014}, where the identity carries a boundary term built on $(u/\delta^s)^2$, with $\delta(x)=\dist(x,\partial\om)$, and where they crucially used the fact that bounded solutions $u$ satisfy the boundary behavior $u/\delta^s\in C^{\alpha}(\bar\om)$ for some $\alpha>0$. On the full space, no boundary term appears, and the identity for the fractional $q$-Laplacian was settled by Ambrosio \cite{pFractionallaplacianPohozaev} for solutions in $W^{s,q}(\RN)\cap C^{1,1}(\RN)$ which solve the equation pointwise, and used it to produce a ground state for a Berestycki-Lions type nonlinearity.
\medskip

Once purely local and nonlocal Pohozaev identities were there, it started the path for the mixed operator $\mathcal L_{2,2}=-\Delta+(-\Delta)^s$. Again Ros-Oton and Serra \cite{RosOtonSerra_PohozaevMixed} proved this identity, on bounded star-shaped domains and for operators which are the sum of a second-order and an integro-differential part. The identity can also be shown to hold for anisotropic nonlocal operators, whose kernels are homogeneous but not radial, and it was done in \cite{RosOtonSerraValdinoci_PohozaevMixedAnisotropic}. With various generalizations and variations on the operator or the nonlinearity, the identity was proved in \cite{Biswas_PohozaevMixed,Anthal_PohozaevChoquard}. Finally, for general $p>1$, this was done by Anthal and Garain \cite{anthal2025pohozaev}, for equations and systems on $\RN$ driven by the anisotropic $p$-Laplacian, by the fractional $p$-Laplacian, and by their sum. In their mixed operator both parts carry the same exponent $p$, and $u$ is assumed H\"older continuous of some order $l>s$. Very recently, Perera and Sportelli \cite{perera2026scaling} considered the sum of two fractional $p$-Laplacians with different orders and different exponents, and proved a Pohozaev identity for finite-energy weak solutions of $(-\Delta)^{s_1}_{p_1}u+(-\Delta)^{s_2}_{p_2}u=\lambda|x|^{-b}|u|^{q-2}u$ without any additional regularity, by a method based on scaling and convexity defects. They also indicate that the argument extends to the mixed case $s_2=1$. Their method relies on the Hardy weight $|x|^{-b}$ with $2b>\max\{q,2\}$, and the unweighted case $b=0$, which is the one considered here, is left open there \cite[Remark~2.8]{perera2026scaling}. For the purely local weighted double phase operator $-\Delta_p u-\dive(|x|^{-b}|\nabla u|^{q-2}\nabla u)$, a Pohozaev identity for all radial finite-energy solutions with weighted power nonlinearities was obtained by Perera and Winkert \cite{perera2026doublephase}, through the radial form of the equation. In our case one part is local and the other nonlocal, and they carry different exponents $p$ and $q$, hence different homogeneities. H\"older continuity is needed only when $p\le sq$, and then only of order $l>\frac{sq-p}{q-p}$, which is smaller than $s$. It would be interesting to know whether the approach of \cite{perera2026scaling} can be adapted to unweighted problems for $\Lpq$ to remove this assumption.

\medskip
The main difficulty in proving the Pohozaev identity is that $x\cdot\nabla u$ is in general not an admissible test function. As in \cite{anthal2025pohozaev,pFractionallaplacianPohozaev,Liu_plapPohozaev}, we test the equation with $\psi_\lambda(x)\sum_jx_jD_ju$, where $D_ju$ is a difference quotient and $\psi_\lambda(x)=\psi(\lambda x)$ is a cut-off function, and we let first $h\to0$ and then $\lambda\to0$. This works because the difference quotient is an admissible test function, and because the convexity of $t\mapsto|t|^p$ and $t\mapsto|t|^q$ gives the remainder terms a sign. The new point in the mixed setting is the term in which the difference quotient of $u$ is paired with the nonlocal kernel. Its treatment depends on the position of $p$ with respect to $q$ and $sq$, so we distinguish the three cases $q<p$, $sq<p\le q$ and $p\le sq$, and the local H\"older continuity of $u$ is used only when $p\le sq$.

\medskip
Once we have the Pohozaev identity, we combine it with the Nehari identity to obtain the expected nonexistence result. Let
\begin{equation*}
    m_-=\min\{p^*,q^*_s\} \quad \text{and} \quad  m_+:=\max\{p^*,q^*_s\}.
\end{equation*}
Nonexistence then holds for nonlinearities which are supercritical for both operators, that is, $tf(t)\ge m_+F(t)$ for all $t$, or subcritical for both operators, that is, $tf(t)\le m_-F(t)$ for all $t$. When $p^*=q^*_s$, these inequalities are required to be strict for $t\neq0$. This leads to the following result.

\begin{theorem}[Nonexistence]\label{thm:nonex}
Assume \eqref{assump: standing} and \eqref{assump:f}. Let $u\in\X(\RN)\cap L^\infty(\RN)\cap C(\RN)$ be a weak solution of \eqref{eq:P} which satisfies \eqref{eq:integrability} and $f(u)u\in L^1(\RN)$. If $p\le sq$, assume also that $u\in C^{0,l}_{\mathrm{loc}}(\RN)$ for some $l\in(l_*,1]$, with $l_*$ as in \eqref{eq:lstar}. Suppose that one of the following holds:
\begin{enumerate}[label=\textup{(\roman*)}]
  \item $t f(t)\ge m_+F(t)$ for all $t\in\R$, and when $\ps=q^*_s$ we require the strict inequality $tf(t)>m_+F(t)$ for all $t\neq0$;
  \item $t f(t)\le m_-F(t)$ for all $t\in\R$, and when $\ps=q^*_s$ we require the strict inequality $tf(t)<m_-F(t)$ for all $t\neq0$.
\end{enumerate}
Then $u\equiv0$.
\end{theorem}
\begin{remark}
(a) If $f$ satisfies \eqref{eq:f-near-zero} for some $a\ge\max\{p,m_-\}$, then the conditions \eqref{eq:integrability} and $f(u)u\in L^1(\RN)$ hold for every $u\in\X(\RN)\cap L^\infty(\RN)$.

(b) Theorem~\ref{thm:nonex} also covers nonlinearities which are not pure powers. For example, let $f(t)=\sum_{k=1}^K a_k|t|^{r_k-2}t$ with $a_k>0$ and $r_k>1$. Then $tf(t)=\sum_k a_k|t|^{r_k}$ and $F(t)=\sum_k \frac{a_k}{r_k}|t|^{r_k}$. Hence (i) holds if $r_k\ge m_+$ for all $k$, and (ii) holds if $r_k\le m_-$ for all $k$ (with $r_k>m_+$, respectively $r_k<m_-$, for all $k$ when $\ps=q^*_s$). Another example is $f(t)=|t|^{r-2}t\,e^{t^2}$ with $r\ge m_+$. Since $f$ is odd, $F$ is even, and for $t\neq0$
\[
F(t)=\int_0^{|t|}\tau^{r-1}e^{\tau^2}\,d\tau<\frac{1}{r}\,|t|^r e^{t^2}=\frac1r\,tf(t),
\]
so (i) holds with strict inequality, also when $\ps=q^*_s$. Here \eqref{eq:f-near-zero} holds with $a=r\ge m_+>p$, so by (a) the integrability conditions are satisfied as well.
\end{remark}

For the pure power $f(t)=|t|^{r-2}t$ with $r>1$, the integrability conditions \eqref{eq:integrability}, together with $f(u)u\in L^1(\RN)$, reduce to
\begin{equation}\label{eq:r-integrability}
|u|^r\in L^1(\RN)\qquad\text{and}\qquad |u|^{r-1}\in L^{p'}(\RN).
\end{equation}
Therefore, Theorem~\ref{thm:nonex} gives the following result.

\begin{corollary}\label{cor:nonex-power}
Assume \eqref{assump: standing} and let $r>1$. Suppose that one of the following holds:
\begin{enumerate}[label=\textup{(\roman*)}]
    \item $\ps\ne q^*_s$ and
\[
  r\le\min\{\ps,q^*_s\}
  \qquad\text{or}\qquad
  r\ge\max\{\ps,q^*_s\};
\]
    \item $\ps= q^*_s\neq r$.
\end{enumerate}
Then \eqref{eq:Pr} has no nontrivial weak solution $u\in\X(\RN)\cap L^\infty(\RN)\cap C(\RN)$ satisfying \eqref{eq:r-integrability} and, in the case $p\le sq$, $u\in C^{0,l}_{\mathrm{loc}}(\RN)$ for some $l\in(l_*,1]$. If $r\ge p$, the condition $|u|^{r-1}\in L^{p'}(\RN)$ in \eqref{eq:r-integrability} can be omitted. If $r\ge\min\{\ps,q^*_s\}$, the condition $|u|^r\in L^1(\RN)$ holds automatically.
\end{corollary}
For supersolutions, the relevant threshold is the Serrin exponent rather than the Sobolev exponents, and for sums of operators including $\Lpq$ it was recently determined by Del Pezzo, Pichucho and Quaas \cite{DelPezzoPichuchoQuaas2026serrin}. The above corollary, together with the existence result of \cite{ourPreviousPaper_isocriticalregimemixedlocalnonlocal}, completes the picture for power nonlinearities in the isocritical case \eqref{eq:isocritical-intro}. When $p^*\neq q^*_s$, however, there remains a window, the open interval $m_-<r<m_+$, in which the Pohozaev identity gives no information. In this range the nonlinearity is supercritical with respect to one of the two operators and subcritical with respect to the other. The subcriticality of $r$ with respect to one operator turns out to be the key property. Indeed, the space $\X^{\mathrm{rad}}(\RN)$ of radial functions with finite energy embeds compactly into $L^r(\RN)$ exactly when $m_-<r<m_+$ (see Proposition~\ref{prop:embedding}). In other words, the exponents excluded by the Pohozaev identity are exactly those for which this compactness fails. In the isocritical case, we have $m_-=m_+=p^*$, so the window is empty. The only exponent not excluded is $r=p^*=q^*_s$, for which compactness is lost, and this is why concentration-compactness was needed in \cite{ourPreviousPaper_isocriticalregimemixedlocalnonlocal}. In the window $m_-<r<m_+$, on the other hand, compactness holds, and we expect to find solutions by the direct method of the calculus of variations. To this end, for $r>1$ we consider the constrained minimization problem
\begin{equation}\label{eq:Sr}
\mathcal S_r:=\inf\Big\{\frac1p[u]_{1,p}^p+\frac1q[u]_{s,q}^q:\ u\in\X(\RN)\cap L^r(\RN),\ \|u\|_r=1\Big\},
\end{equation}
and we denote by $\mathcal S$ and $\mathcal S_f$ the best constants in the Sobolev inequalities
\[
\mathcal S\,\|u\|_{p^*}^p\le[u]_{1,p}^p,\qquad \mathcal S_f\,\|u\|_{q^*_s}^q\le[u]_{s,q}^q,\qquad u\in\X(\RN).
\]
We obtain the following result.

\begin{theorem}[existence]\label{thm:existenceNonexistence}
Assume \eqref{assump: standing} and $p^*\neq q^*_s$, let $r>1$, and let $\mathcal S_r$ be defined by \eqref{eq:Sr}.
\begin{enumerate}[label=\textup{(\alph*)}]
\item If $r<m_-$ or $r>m_+$, then $\mathcal S_r=0$, and $\mathcal S_r$ is not attained.
\item If $r=p^*$, then $\mathcal S_r=\frac1p\mathcal S$, and if $r=q^*_s$, then $\mathcal S_r=\frac1q\mathcal S_f$. In both cases $\mathcal S_r$ is not attained.
\item If $m_-<r<m_+$, then $\mathcal S_r>0$ and $\mathcal S_r$ is attained by a nonnegative, radially symmetric and radially nonincreasing function $u$. Moreover, there exists $\mu>0$ such that $u$ satisfies 
\begin{equation}\label{eq:Pmu}
\langle u,\varphi\rangle_p+\langle u,\varphi\rangle_{s,q}=\mu\int_{\RN}|u|^{r-2}u\,\varphi\,dx\qquad\text{for every }\varphi\in\X(\RN).
\end{equation}
Let
\begin{equation}\label{eq:Dr}
    D(r) := r(p-sq) - pq(1-s),
\end{equation}
which is nonzero by Lemma~\ref{lem:Dnonzero}. Then the function $v$ given by 
\begin{equation}\label{eq:normalisation}
v(x):=\kappa\,u(\lam x),\qquad \text{with}\quad \kappa:=\mu^{\frac{p-sq}{D(r)}},\qquad \lam:=\mu^{\frac{q-p}{D(r)}},
\end{equation}
is a nontrivial, nonnegative, radially symmetric and radially nonincreasing weak solution of \eqref{eq:Pr} in $\X(\RN)$.
\end{enumerate}
\end{theorem}
Corollary~\ref{cor:nonex-power} already rules out solutions of \eqref{eq:Pr} for $r\le m_-$ and $r\ge m_+$ within the class considered there, and parts (a) and (b) of Theorem~\ref{thm:existenceNonexistence} explain the variational mechanism behind this. For $r<m_-$ or $r>m_+$ the infimum is zero, so minimizing sequences lose all their energy. At each endpoint, the contribution of one of the two operators vanishes along suitably rescaled sequences, and $\mathcal S_r$ reduces to the sharp Sobolev constant of the other operator. Thus, compactness is lost by concentration at $r=m_+$ and by spreading at $r=m_-$. Inside $(m_-,m_+)$, concentration is prevented by one operator and spreading by the other, and this is why the minimization over the unit sphere of $L^r(\RN)$ succeeds. Finally, $D(r)\neq0$ on $[m_-,m_+]$ when $p^*\neq q^*_s$ (Lemma~\ref{lem:Dnonzero}), and this is what allows the rescaling \eqref{eq:normalisation} to remove the Lagrange multiplier (Lemma~\ref{lem:normalisation}).

\begin{remark}
\begin{enumerate}[label=(\roman*)]
    \item Theorem~\ref{thm:existenceNonexistence}(b) applies to every $u\in\X(\RN)$ without additional regularity assumptions. By itself, non-attainment does not rule out weak solutions of \eqref{eq:Pr} at $r=m_\pm$, since a solution need not always be a minimizer.

    \item When $p=q=2$, we have $2^*>\frac{2N}{N-2s}$, so $r=2^*$ is the upper endpoint. Theorem~\ref{thm:existenceNonexistence}(b) then shows that the sharp Sobolev constant of $-\Delta+(-\Delta)^s$ on $\RN$ coincides with that of $-\Delta$ and is not attained, recovering the known result of \cite{biagiBrezisNirenbergMixed}.
\end{enumerate}
\end{remark}

Next, we look for ground states. For $r>1$, the energy functional and the Pohozaev functional associated with \eqref{eq:Pr} are
\begin{equation}\label{eq:JG-intro}
\mathcal J(u):=\frac1p[u]_{1,p}^p+\frac1q[u]_{s,q}^q-\frac1r\|u\|_r^r,
\qquad
\mathcal G(u):=\frac{N-p}{p}[u]_{1,p}^p+\frac{N-sq}{q}[u]_{s,q}^q-\frac Nr\|u\|_r^r ,
\end{equation}
and we set
\[
\mathcal M:=\{u\in\X(\RN)\setminus\{0\}:\ \mathcal G(u)=0\},\qquad c_{\mathcal M}:=\inf_{\mathcal M}\mathcal J .
\]
For $m_-\le r\le m_+$ these functionals are well defined and of class $C^1$ on $\X(\RN)$ by Lemma~\ref{lem:interp}. A nontrivial weak solution of \eqref{eq:Pr} lies in the Pohozaev set $\mathcal M$ if and only if it satisfies the Pohozaev identity \eqref{eq:poho}. By a \emph{ground state} we mean a nontrivial weak solution $w$ of \eqref{eq:Pr} such that $w\in\mathcal M$ and $\mathcal J(w)=c_{\mathcal M}$. Thus, a ground state has the least energy among all nontrivial solutions which satisfy the Pohozaev identity.

\medskip
The solution in Theorem~\ref{thm:existenceNonexistence} is obtained by minimizing a quotient, so it is not clear a priori whether it is a ground state. To construct a ground state, we use $\mathcal{M}$ as a natural constraint. Since $\mathcal J-\frac1N\mathcal G=\frac1N\big([u]_{1,p}^p+s[u]_{s,q}^q\big)$, the energy $\mathcal{J}$ is positive and coercive on $\mathcal{M}$. On the Nehari manifold, by contrast, this can fail when $r<\max\{p,q\}$, a case which does occur within $(m_-,m_+)$. In that case $\mathcal J(tu)\to+\infty$ as $t\to\infty$, so the functional has no mountain-pass geometry along the rays $t\mapsto tu$. To bypass this difficulty, we follow the approach of Jeanjean \cite{jeanjean_1997} (see also \cite{hirata_ikoma_tanaka_2010}) and work with the dilations $u_t(x):=u(x/t)$. Since $[u_t]_{1,p}^p=t^{N-p}[u]_{1,p}^p$, $[u_t]_{s,q}^q=t^{N-sq}[u]_{s,q}^q$ and $\|u_t\|_r^r=t^N\|u\|_r^r$, the $L^r$ term always dominates along dilations for large $t$. This allows us to construct a Palais--Smale sequence which asymptotically satisfies the Pohozaev identity (Lemma~\ref{lem:jeanjean}). For the purely local operator $-\Delta_p-\Delta_q$, a similar combination of a Pohozaev identity and Jeanjean's argument was used by Pomponio and Watanabe \cite{pomponio_watanabe_2018} for zero-mass problems (see also \cite{demedeiros_nascimento_siciliano_2026}).
\begin{theorem}\label{thm:groundstate}
Assume \eqref{assump: standing}, $\ps\neq q^*_s$ and $m_-<r<m_+$. Then $c_{\mathcal M}>0$ is attained by a nonnegative, radially symmetric weak solution $w\in\X(\RN)$ of \eqref{eq:Pr}, that is, $w$ is a ground state: $\inf_{\mathcal M}\mathcal J=\mathcal J(w)$. Moreover, $c_{\mathcal M}=c_{\mathrm{mp}}$, where $c_{\mathrm{mp}}$ is the mountain-pass level defined in \eqref{eq:cmp}. In particular, $\mathcal J(w)\le\mathcal J(z)$ for every nontrivial weak solution $z\in\X(\RN)$ of \eqref{eq:Pr} which satisfies the Pohozaev identity \eqref{eq:poho}.
\end{theorem}
\begin{remark}[Positivity]
If $p\le q$, the solutions $v$ of Theorem~\ref{thm:existenceNonexistence} and $w$ of Theorem~\ref{thm:groundstate} are positive a.e.\ in $\RN$. Indeed, they are nontrivial, nonnegative weak supersolutions of $\Lpq$ in $\RN$, because the right-hand side of \eqref{eq:Pr} is nonnegative for nonnegative $u$. Hence the strong maximum principle of \cite{ourPreviousPaper_isocriticalregimemixedlocalnonlocal} applies, and it holds for $1<p\le q$ without any relation between $\ps$ and $q^*_s$. 
\end{remark}

Theorems~\ref{thm:existenceNonexistence}(c) and~\ref{thm:groundstate} produce solutions of \eqref{eq:Pr} by two different minimizations, and it is natural to ask how they are related. It turns out that the solution $v$ of Theorem~\ref{thm:existenceNonexistence} is also a ground state. The key point is that the dilations $u^{(\eta)}:=\eta^{N/r}u(\eta\,\cdot\,)$, $\eta>0$, preserve the constraint $\|u\|_r=1$. Hence, for the minimizer $u$ of Theorem~\ref{thm:existenceNonexistence}, the map $\eta\mapsto\frac1p[u^{(\eta)}]_{1,p}^p+\frac1q[u^{(\eta)}]_{s,q}^q$ attains its minimum at $\eta=1$. This gives the Pohozaev identity for the equation \eqref{eq:Pmu} satisfied by $u$, and therefore $v\in\mathcal M$, without any regularity assumption. Moreover, both minimizations are governed by the best constant in an interpolation inequality. For $m_-<r<m_+$, let $\theta\in(0,1)$ be defined by $\frac1r=\frac{\theta}{p^*}+\frac{1-\theta}{q^*_s}$, and let $\mathcal C_r$ be the largest constant such that
\[
\mathcal C_r\,\|u\|_r\le[u]_{1,p}^{\theta}\,[u]_{s,q}^{1-\theta}\qquad\text{for all }u\in\X(\RN);
\]
which is positive by the H\"older and Sobolev inequalities (Lemma~\ref{lem:interp}). Finally, set
\[
\frac{1}{p_\theta}:=\frac{\theta}{p}+\frac{1-\theta}{q},\qquad c_\theta:=\theta^{\theta/p}(1-\theta)^{(1-\theta)/q} .
\]
We have the following result.

\begin{theorem}\label{thm:same}
Assume \eqref{assump: standing}, $\ps\neq q^*_s$ and $m_-<r<m_+$, and let $v$ be the solution of Theorem~\ref{thm:existenceNonexistence}. Then $\mathcal J(v)=c_{\mathcal M}$, i.e., $v$ is a nonnegative, radially nonincreasing ground state of \eqref{eq:Pr}. Moreover, the Lagrange multiplier in Theorem~\ref{thm:existenceNonexistence}
is $\mu=p_\theta\,\mathcal S_r$, and
\begin{equation}\label{eq:cM-Sr}
c_{\mathcal M}
=\frac{r-p_\theta}{r\,p_\theta}\,\big(p_\theta\,\mathcal S_r\big)^{\frac{r}{r-p_\theta}}
=\frac{r-p_\theta}{r\,p_\theta}\Big(\frac{\mathcal C_r}{c_\theta}\Big)^{\frac{r\,p_\theta}{r-p_\theta}} .
\end{equation}
\end{theorem}
The paper is organized as follows. Section~\ref{sec:prelim} contains notation and the function spaces. In Section~\ref{sec:pohozaev}, we prove the Pohozaev identity, and in Section~\ref{sec:nonex} we deduce the nonexistence results. Finally, in Section~\ref{sec:existence} we prove the existence results in the range $m_-<r<m_+$.

\section{Notations and preliminaries}\label{sec:prelim}

\subsection{Notations}

Throughout the paper, we use the following notations.
\begin{itemize}
\item $B_R(x_0)=\{x\in\RN:|x-x_0|<R\}$, $B_R:=B_R(0)$. 
\item For a set $A \subset \mathbb{R}^N$, $
A^c := \mathbb{R}^N \setminus A.
$
\item For any measurable set $E \subset \mathbb{R}^N$, $|E|$ denotes its Lebesgue measure and $\mathbf 1_E$ its indicator.
\item $\mathbb S^{N-1}$ is the unit sphere and $|\mathbb S^{N-1}|:=\mathcal H^{N-1}(\mathbb S^{N-1})$,
  where $\mathcal H^{N-1}$ is the $(N-1)$-dimensional Hausdorff measure.
\item For a measurable set $\Omega\subseteq\RN$ and $1\le t\le\infty$, $L^t(\Omega)$ denotes the usual
  Lebesgue space with norm $\|\cdot\|_{L^t(\Omega)}$, and $\|\cdot\|_t:=\|\cdot\|_{L^t(\RN)}$.
  For $1<t<\infty$, $t':=\frac{t}{t-1}$.
\item For $\Omega$ open, $W^{1,t}_{\rm loc}(\Omega)$ is the usual local Sobolev space. $\nabla u$ and
  $u_{x_j}=\partial_{x_j}u$ always denote weak derivatives, and $D^2u$ the Hessian.
  $C^k_c(\Omega)$ denotes $C^k$ functions with compact support in $\Omega$, and
  $\|g\|_{C^2(K)}:=\sum_{|\beta|\le2}\|\partial^\beta g\|_{L^\infty(K)}$.
\item For $l\in(0,1]$ and $K\subset\RN$ compact, $ [u]_{C^{0,l}(K)}:=\sup_{\substack{x,y\in K\\ x\neq y}}\frac{|u(x)-u(y)|}{|x-y|^{l}},$
  \[
  C^{0,l}_{\rm loc}(\RN):=\big\{u\in C(\RN):\ [u]_{C^{0,l}(K)}<\infty\ \text{for every compact }K\subset\RN\big\}.
  \]
\item For a measurable function $u$, $u^{\pm}:=\max\{\pm u,0\}$, so that $u=u^+-u^-$.
\item For $u:\RN\to\R$ measurable,
  \[
  [u]_{1,p}:=\Big(\int_{\RN}|\nabla u|^p\,dx\Big)^{1/p},\qquad
  [u]_{s,q}:=\Big(\iint_{\R^{2N}}|u(x)-u(y)|^q\,d\mu_{s,q}\Big)^{1/q},
  \]
  where $d\mu_{s,q}:=|x-y|^{-N-sq}\,dx\,dy$.
\item For $t>1$, we define $J_t:\mathbb{R}\to\mathbb{R}$ by $J_t(\tau):=|\tau|^{t-2}\tau$.
\item $C$ denotes a positive constant depending only on $N,p,q,s$, whose value may change from line
  to line. Moreover, the dependence on further quantities is indicated, e.g.\ $C(\lambda)$.
  For $a,b>0$, $a\gg b$ means $a\ge Cb$ for a sufficiently large such constant, for instance, $k\gg1$ means that $k$ is sufficiently large.
\end{itemize}

\subsection{The operator and the function spaces}
For $1<p<N$, $D^{1,p}_0(\RN)$ is the completion of $C^\infty_c(\RN)$ with respect to $[\cdot]_{1,p}$, and
for $s\in(0,1)$, $sq<N$, $D^{s,q}_0(\RN)$ is the completion with respect to $[\cdot]_{s,q}$ (for more details on the spaces involved refer to \cite{ourPreviousPaper_isocriticalregimemixedlocalnonlocal}). We set
\[
\X(\RN):=D^{1,p}_0(\RN)\cap D^{s,q}_0(\RN),\qquad \|u\|_{\X}:=[u]_{1,p}+[u]_{s,q},
\]
By the Sobolev inequalities, $\X(\RN)\hookrightarrow L^{p^*}(\RN)\cap L^{q^*_s}(\RN)$
(see \cite{DNPV2012} for the fractional case).
For $u,\varphi\in\X(\RN)$, we set
\[
\langle u,\varphi\rangle_{p}:=\int_{\RN}|\nabla u|^{p-2}\nabla u\cdot\nabla\varphi\,dx,
\qquad
\langle u,\varphi\rangle_{s,q}:=\iint_{\R^{2N}}J_q\big(u(x)-u(y)\big)\big(\varphi(x)-\varphi(y)\big)\,d\mu_{s,q}.
\]
By H\"older's inequality,
$|\langle u,\varphi\rangle_{p}|\le[u]_{1,p}^{p-1}[\varphi]_{1,p}$ and
$|\langle u,\varphi\rangle_{s,q}|\le[u]_{s,q}^{q-1}[\varphi]_{s,q}$, thus both are well defined,
and linear and continuous in $\varphi$. Accordingly, the operators
$-\Delta_p,\,(-\Delta)^s_q:\X(\RN)\to\X(\RN)^*$ are defined by
\[
\langle-\Delta_pu,\varphi\rangle:=\langle u,\varphi\rangle_{p},\qquad
\langle(-\Delta)^s_qu,\varphi\rangle:=\langle u,\varphi\rangle_{s,q},
\]
where $\langle\cdot,\cdot\rangle$ is the duality pairing between $\X(\RN)^*$ and $\X(\RN)$, and we write
\[
\Lpq:=-\Delta_p+(-\Delta)^s_q,\qquad
\langle\Lpq u,\varphi\rangle=\langle u,\varphi\rangle_{p}+\langle u,\varphi\rangle_{s,q}.
\]

\section{The Pohozaev identity}\label{sec:pohozaev}

In this section, we prove Theorem~\ref{thm:main Pohozaev}, that is, the Pohozaev identity for weak solutions of \eqref{eq:P}, where the nonlinearity $f$ satisfies \eqref{assump:f}.
 
\subsection{Difference quotients}
Throughout, $\nabla u=(u_{x_1},\dots,u_{x_N})$ and each $u_{x_j}=\partial_{x_j}u$ denote the weak (distributional) partial derivatives: for $u\in D^{1,p}_0(\RN)$, $\partial_{x_j}u$ is the unique $L^p(\RN)$ function satisfying
\[
\int_{\RN}u\,\partial_{x_j}\varphi\,dx=-\int_{\RN}(\partial_{x_j}u)\,\varphi\,dx
\]
for all $\varphi\in C_c^\infty(\RN)$.
For $v:\RN\to\R$ and $h\ne0$, the (forward) difference quotient in the $j$-th direction is
\[
    D_jv(x)=\frac{v(x+he_j)-v(x)}{h},
    \qquad
    \delta_jv(x)=v(x+he_j),
\]
where $e_j$ is the $j$-th coordinate unit vector. The shift operator $\delta_j$ satisfies the discrete Leibniz rule
\begin{equation*}
    D_j(vw)=w\,D_jv+(\delta_jv)\,D_jw,
    \qquad j=1,\dots,N. 
\end{equation*}
In the proof of Theorem~\ref{thm:main Pohozaev}, we always take $h>0$, with $h\to0$ understood as $h\to0^+$.

\begin{lemma}\label{lem:null Identity}
Let $h\ne0$ and $j\in\{1,\dots,N\}$, and write $\delta_j U(x)=U(x+he_j)$.
\begin{enumerate}
\item[\textup{(a)}] For every $U\in L^1(\RN)$,
    \begin{equation*}
        \int_{\RN}D_jU\,dx=0 .
    \end{equation*}
    For every $G\in L^1(\R^{2N},\dmu)$,
    \begin{equation*}
        \iint_{\R^{2N}}D_j^{\mathrm{diag}}G\,\dmu=0,
        \qquad
        D_j^{\mathrm{diag}}G(x,y):=\frac{G(x+he_j,y+he_j)-G(x,y)}{h}.
    \end{equation*}
\item[\textup{(b)}] For $U\in L^1(\RN)$ and $g\in L^\infty(\RN)$,
    \begin{equation}\label{eq:dibp-loc}
        \int_{\RN}g\,D_jU\,dx=-\int_{\RN}(\delta_jU)\,D_jg\,dx .
    \end{equation}
    For $G\in L^1(\R^{2N},\dmu)$,
    \begin{equation}\label{eq:dibp-nl}
        \iint_{\R^{2N}}g(x)\,D_j^{\mathrm{diag}}G\,\dmu
        =-\iint_{\R^{2N}}(\delta_jG)\,D_jg(x)\,\dmu,
        \quad \delta_jG(x,y):=G(x+he_j,y+he_j).
    \end{equation}
\end{enumerate}
\end{lemma}
\begin{proof}
(a) Follows from the translation invariance of the Lebesgue measure. For the second identity, note that the measure $\dmu$ is invariant under the diagonal translation $(x,y)\mapsto(x+he_j,y+he_j)$.\\
(b) Using the same invariance in the term containing $g(x)U(x)$, we get
\[
\begin{aligned}
\int_{\RN}g\,D_jU\,dx&=\frac1h\int_{\RN}g(x)U(x+he_j)\,dx-\frac1h\int_{\RN}g(x+he_j)U(x+he_j)\,dx\\
&=-\int_{\RN}(\delta_jU)\,D_jg\,dx,
\end{aligned}
\]
and \eqref{eq:dibp-nl} is proved in the same way. All integrals are finite since $g\in L^\infty$, $U\in L^1(\RN)$, and $G\in L^1(\R^{2N},\dmu)$.
\end{proof}

\subsection{The Pohozaev identity}

Let $f$ and $F$ be as in \eqref{assump:f}, and let $u$ be as in Theorem~\ref{thm:main Pohozaev}.
Fix $\psi\in C^1_c(\RN)$ with $0\le\psi\le1$, $\psi\equiv1$ on $\{|x|\le1\}$,
and $\psi\equiv0$ on $\{|x|\ge2\}$. For $\lam>0$, set $\psi_\lam(x)=\psi(\lam x)$.
Then, with a constant $C>0$ independent of $\lam$,
\begin{equation}\label{eq:cutoff}
    0\le\psi_\lam\le1,
    \qquad
    |x|\,|\nabla\psi_\lam(x)|\le C\quad\text{for all }x\in\RN,\ \lam>0.
\end{equation}
Define the Lipschitz, compactly supported vector field
\begin{equation*}
    \Psi_\lam(x)=\psi_\lam(x)\,x.
\end{equation*}
From \eqref{eq:cutoff} one obtains, for a constant $C>0$ independent of $\lam$,
\[
    |\Psi_\lam(x)-\Psi_\lam(y)|\le C|x-y|,
    \qquad
    \Psi_\lam\in C^{0,1}_c(\RN,\RN).
\]
The Pohozaev test function is 
\[
    \varphi
    =\psi_\lam\sum_{j=1}^N x_j\,D_j u
    \;=\;\sum_{j=1}^N \psi_\lam\,x_j\,D_j u .
\]
Since $u\in\X(\RN)\cap L^\infty(\RN)$, each $D_ju$ belongs to $\X(\RN)\cap L^\infty(\RN)$, and $\psi_\lam x_j$ is a bounded Lipschitz function. Consequently, $\varphi\in \X(\RN)$. Moreover, $\varphi$ has compact support, and since $f(u)\in L^\infty(\RN)$, we have $f(u)\varphi\in L^1(\RN)$. Hence $\varphi$ is admissible in \eqref{eq:weak}, so
\begin{equation}\label{eq:master}
    \underbrace{\int_{\RN}|\nabla u|^{p-2}\nabla u\cdot\nabla\varphi\,dx}_{=:A}
    +\underbrace{\iint_{\R^{2N}}\Jq(u(x)-u(y))(\varphi(x)-\varphi(y))\dmu}_{=:B}
    =\underbrace{\int_{\RN}f(u)\varphi\,dx}_{=:L}.
\end{equation}
We treat $A$ (local), $B$ (nonlocal), and $L$ (right-hand side) separately, letting first $h\to0$ and subsequently $\lam\to0$.

\subsubsection{The local term A}
We expand $A=A_1+A_2+A_3$ as
\begin{align}
  A_1&:=\sum_{i,j=1}^N\int_{\RN}|\nabla u|^{p-2}u_{x_i}\,
        (\partial_{x_i}\psi_\lam)\,x_j\,D_ju\,dx, \label{eq:A1}\\
  A_2&:=\sum_{i,j=1}^N\int_{\RN}|\nabla u|^{p-2}u_{x_i}\,
        \psi_\lam\,\delta_{ij}\,D_ju\,dx, \notag\\
  A_3&:=\sum_{j=1}^N\int_{\RN}|\nabla u|^{p-2}\,\nabla u\cdot D_j(\nabla u)\,
       \psi_\lam x_j\,dx=A_{3,1}-A_{3,2}, \notag
\end{align}
where
\begin{align*}
A_{3,1}&:=\frac1p\sum_{j=1}^N\int_{\RN}D_j\big(|\nabla u|^p\big)\psi_\lam x_j\,dx, \\
  A_{3,2}&:=\frac1p\sum_{j=1}^N\int_{\RN}
      \Big(D_j(|\nabla u|^p)-p\,|\nabla u|^{p-2}\nabla u\cdot D_j(\nabla u)\Big)
      \psi_\lam x_j\,dx. 
\end{align*}
 
\paragraph{Limit of $A_1$.}
Recall from \eqref{eq:A1}
\begin{equation*}
  A_1=A_1(h)=\sum_{i,j=1}^N\int_{\RN}
     |\nabla u(x)|^{p-2}u_{x_i}(x)\,\big(\partial_{x_i}\psi_\lam(x)\big)\,x_j\,
     D_ju(x)\,dx,
\end{equation*}
and that the only factor depending on $h$ is the difference quotient
$D_ju(x)=\dfrac{u(x+he_j)-u(x)}{h}$. We prove
\begin{equation*}
  \lim_{h\to0}A_1(h)
  =\int_{\RN}|\nabla u|^{p-2}(\nabla u\cdot\nabla\psi_\lam)(x\cdot\nabla u)\,dx .
\end{equation*}
 
\smallskip
For fixed $i,j\in\{1,\dots,N\}$ define
\begin{equation*}
  G_{ij}(x):=|\nabla u(x)|^{p-2}\,u_{x_i}(x)\,\big(\partial_{x_i}\psi_\lam(x)\big)\,x_j,
\end{equation*}
so that $A_1(h)=\sum_{i,j}\int_{\RN}G_{ij}(x)\,D_ju(x)\,dx$. The candidate limit
is $A_1^{0}:=\sum_{i,j}\int_{\RN}G_{ij}(x)\,u_{x_j}(x)\,dx$, and we must show
$A_1(h)\to A_1^0$. Observe that $G_{ij}\in L^{p'}(\RN)$, where $p'=\tfrac{p}{p-1}$. Therefore, it suffices to show that 
\begin{equation}\label{eq:Djconv}
  \big\|D_ju-u_{x_j}\big\|_{p}\xrightarrow[h\to0]{}0 ,
  \qquad D_ju=\frac{u(\cdot+he_j)-u(\cdot)}{h}.
\end{equation}

Since $u\in W^{1,p}_{\mathrm{loc}}(\RN)$ with $\partial_{x_j}u\in L^p(\RN)$, by \cite[Theorem $2.36$]{kinnunen2017sobolev}, for every fixed $h\ne0$,
\begin{equation}\label{eq:rep}
  D_ju(x)=\frac{u(x+he_j)-u(x)}{h}
        =\int_0^1\partial_{x_j}u(x+the_j)\,dt
  \qquad\text{for a.e.\ }x\in\RN .
\end{equation}
 
Further, recall that for any $g\in L^p(\RN)$, $\|g(\cdot+\xi)-g\|_{p}\to0$ as
$\xi\to0$. 
Apply \eqref{eq:rep} 
\begin{equation*}
  D_ju(x)-u_{x_j}(x)
  =\int_0^1\big(u_{x_j}(x+the_j)-u_{x_j}(x)\big)\,dt .
\end{equation*}
By Minkowski's integral inequality, $\|D_ju-u_{x_j}\|_{p}\le\int_0^1\|u_{x_j}(\cdot+the_j)-u_{x_j}\|_{p}\,dt$, so by the
dominated convergence theorem
\[
  \int_0^1\big\|u_{x_j}(\cdot+the_j)-u_{x_j}\big\|_{p}\,dt\xrightarrow[h\to0]{}0 ,
\]
which proves \eqref{eq:Djconv}.
Hence
\[
  \int_{\RN}G_{ij}\,D_ju\,dx\xrightarrow[h\to0]{}\int_{\RN}G_{ij}\,u_{x_j}\,dx
  \quad\text{for every fixed }i,j .
\]
Therefore
\begin{equation*}
  \lim_{h\to0}A_1
  =\int_{\RN}|\nabla u|^{p-2}\big(\nabla u\cdot\nabla\psi_\lam\big)
       \big(x\cdot\nabla u\big)\,dx.
\end{equation*}
\paragraph{Limit of $A_2$.}
\[
  A_2(h)=\sum_{j=1}^N\int_{\RN}\widehat H_j(x)\,D_ju(x)\,dx,
  \qquad \widehat H_j:=|\nabla u|^{p-2}u_{x_j}\,\psi_\lam.
\]
Applying H\"older with
exponents $p',p$, and \eqref{eq:Djconv} give
\begin{equation*}
  \lim_{h\to0}A_2
  =\sum_{j=1}^N\int_{\RN}|\nabla u|^{p-2}u_{x_j}\,\psi_\lam\,u_{x_j}\,dx
  =\int_{\RN}\psi_\lam|\nabla u|^{p}\,dx.
\end{equation*}
 
\paragraph{Limit of $A_{3,1}$.}
By the discrete integration by parts \eqref{eq:dibp-loc}, applied with $U=|\nabla u|^p\in L^1(\RN)$ and $g=\psi_\lam x_j$,
\[
A_{3,1}=-\frac1p\sum_{j=1}^N\int_{\RN}|\nabla u(x+he_j)|^p\,D_j(\psi_\lam x_j)\,dx .
\]
Since $|\nabla u(\cdot+he_j)|^p\to|\nabla u|^p$ in $L^1(\RN)$, while $D_j(\psi_\lam x_j)\to\partial_{x_j}(\psi_\lam x_j)$ uniformly on $\RN$ and $\|D_j(\psi_\lam x_j)\|_{L^\infty}$ is bounded independently of $h$, we obtain
\begin{align*}
  \lim_{h\to0}A_{3,1}
  =-\frac1p\sum_{j=1}^N\int_{\RN}|\nabla u|^p\,
       \partial_{x_j}(\psi_\lam x_j)\,dx &=-\frac1p\int_{\RN}|\nabla u|^p\Big(\underbrace{\textstyle\sum_j\partial_{x_j}(\psi_\lam)x_j}_{\nabla\psi_\lam\cdot x}+\underbrace{\textstyle\sum_j\psi_\lam}_{N\psi_\lam}\Big)dx \notag\\
  &=-\frac{N}{p}\int_{\RN}|\nabla u|^p\psi_\lam\,dx
    -\frac1p\int_{\RN}|\nabla u|^p\,\nabla\psi_\lam\cdot x\,dx.
\end{align*}
\subsubsection{The nonlocal term B}
\begin{align*}
  B&:=\iint_{\R^{2N}}\Jq(u(x)-u(y))\big(\varphi(x)-\varphi(y)\big)\dmu\\
  & =\sum_{j=1}^N\iint_{\R^{2N}}\Jq(u(x)-u(y))
     \big(\psi_\lam(x)x_j D_ju(x)-\psi_\lam(y)y_j D_ju(y)\big)\dmu :=B_1+B_2,\notag
\end{align*}
where
\begin{align*}
  B_1&:=\sum_{j=1}^N\iint_{\R^{2N}}\Jq(u(x)-u(y))
        \big(\psi_\lam(x)x_j-\psi_\lam(y)y_j\big)D_ju(x)\dmu, \\
  B_2&:=\sum_{j=1}^N\iint_{\R^{2N}}\Jq(u(x)-u(y))
        \big(D_ju(x)-D_ju(y)\big)\psi_\lam(x)x_j\dmu .
\end{align*}
Further,
$B_2=B_{2,1}-B_{2,2}$ with
\begin{align*}
  B_{2,1}&:=\frac1q\sum_{j=1}^N\iint_{\R^{2N}}
     D_j^{\mathrm{diag}}\big(|u(x)-u(y)|^q\big)\,\psi_\lam(x)x_j\dmu,\\
  B_{2,2}&:=\frac1q\sum_{j=1}^N\iint_{\R^{2N}}
     R^{\mathrm{nl}}_j(x,y)\,\psi_\lam(x)x_j\dmu,
\end{align*}
where $R^{\mathrm{nl}}_j:=D_j^{\mathrm{diag}}\big(|u(x)-u(y)|^q\big)-q\,\Jq(u(x)-u(y))\big(D_ju(x)-D_ju(y)\big)$. Thus $B=B_1+B_{2,1}-B_{2,2}.$

\smallskip
\textbf{Claim. }$\lim_{h\to0}(A_{3,2}+B_{2,2})=0$.

\smallskip
Recall
\begin{align*}
A_{3,2}=&\frac1p\sum_{j=1}^N\int_{\RN}R^{\mathrm{loc}}_j(x)\,\psi_\lam(x)x_j\,dx\,,
  \quad R^{\mathrm{loc}}_j:=D_j(|\nabla u|^p)-p\,|\nabla u|^{p-2}\nabla u\cdot D_j(\nabla u).
\end{align*}
By convexity of $t\mapsto|t|^p$ and $t\mapsto|t|^q$, for $h>0$ and
every $x,y$,$\,\,R^{\mathrm{loc}}_j(x)\ge0,
  \,\, R^{\mathrm{nl}}_j(x,y)\ge0.$
Since the integrand remainders are
nonnegative,
\begin{equation}\label{eq:pullweight}
  \big|A_{3,2}+B_{2,2}\big|
  \le C(\lam)\sum_{j=1}^N\Big[\frac1p\int_{\RN}R^{\mathrm{loc}}_j\,dx
       +\frac1q\iint_{\R^{2N}}R^{\mathrm{nl}}_j\dmu\Big].
\end{equation}
Using the identities of Lemma~\ref{lem:null Identity}\,(a), we get 
\begin{align*}
  &\frac1p\int R^{\mathrm{loc}}_j\,dx +\frac1q\iint R^{\mathrm{nl}}_j\dmu \\
  &=-\Big[\int_{\RN}|\nabla u|^{p-2}\nabla u\cdot\nabla(D_ju)\,dx
   +\iint_{\R^{2N}}\Jq(u(x)-u(y))\big(D_ju(x)-D_ju(y)\big)\dmu\Big] \\
  &=-\int_{\RN}f(u)\,D_ju\,dx.
\end{align*}
In the last step, we used \eqref{eq:weak} with $\varphi=D_ju$. This choice is admissible because $D_ju\in\X(\RN)$, and $f(u)D_ju\in L^1(\RN)$ since $f(u)\in L^{p'}(\RN)$ by \eqref{eq:integrability} and $D_ju\in L^p(\RN)$. Indeed, by \eqref{eq:rep}, Minkowski's integral inequality and translation invariance,
\begin{equation*}
  \|D_ju\|_{L^p(\RN)}
  =\Big\|\int_0^1\partial_{x_j}u(\cdot+the_j)\,dt\Big\|_{L^p(\RN)}
  \le\int_0^1\|\partial_{x_j}u(\cdot+the_j)\|_{L^p(\RN)}\,dt
  =\|\partial_{x_j}u\|_{L^p(\RN)} .
\end{equation*}
Inserting the identity above into \eqref{eq:pullweight}, we obtain
\begin{equation}\label{eq:A32lim}
  \big|A_{3,2}+B_{2,2}\big|
  \le C(\lam)\sum_{j=1}^N\Big|\int_{\RN}f(u)\,D_ju\,dx\Big| .
\end{equation}
Now let $h\to0$. Since $f(u)\in L^{p'}(\RN)$ and $D_ju\to u_{x_j}$ in $L^p(\RN)$ by \eqref{eq:Djconv}, H\"older's inequality gives
\[
  \lim_{h\to0}\int_{\RN}f(u)\,D_ju\,dx=\int_{\RN}f(u)\,u_{x_j}\,dx .
\]
We show that the last integral is zero. Let $M:=\|u\|_\infty$. Since $F\in C^1(\R)$ and $F'=f$ is bounded on $[-M,M]$, the chain rule for Sobolev functions gives $F(u)\in W^{1,1}_{\mathrm{loc}}(\RN)$ and
\begin{equation}\label{eq:chain}
  \partial_{x_j}\big(F(u)\big)=f(u)\,u_{x_j}\qquad\text{a.e. in }\RN .
\end{equation}
By H\"older's inequality, $f(u)u_{x_j}\in L^1(\RN)$, and $F(u)\in L^1(\RN)$ by \eqref{eq:integrability}. Fix $\chi\in C^\infty_c(\RN)$ with $0\le\chi\le1$ and $\chi\equiv1$ on $B_1$, and set $\chi_R(x):=\chi(x/R)$ for $R>0$. By the definition of the weak derivative,
\[
  \Big|\int_{\RN}\partial_{x_j}\big(F(u)\big)\,\chi_R\,dx\Big|
  =\Big|\int_{\RN}F(u)\,\partial_{x_j}\chi_R\,dx\Big|
  \le\frac{\|\nabla\chi\|_\infty}{R}\,\|F(u)\|_{L^1(\RN)}
  \xrightarrow[R\to\infty]{}0,
\]
while the left-hand side tends to $\big|\int_{\RN}\partial_{x_j}(F(u))\,dx\big|$ by the dominated convergence theorem. Hence
\begin{equation*}
\int_{\RN}f(u)\,u_{x_j}\,dx=\int_{\RN}\partial_{x_j}\big(F(u)\big)\,dx=0 .
\end{equation*}
Therefore, from \eqref{eq:A32lim},
\begin{equation*}
  \lim_{h\to0}\big(A_{3,2}+B_{2,2}\big)=0 .
\end{equation*}
Note that this step uses only \eqref{eq:integrability}, and no growth condition on $f$.

\paragraph{Limit of $B_1$}
\[
  B_1=\sum_{j=1}^N\iint_{\R^{2N}}\Jq(u(x)-u(y))\,
       w_j(x,y)\,D_ju(x)\dmu,
  \quad w_j(x,y):=\psi_\lam(x)x_j-\psi_\lam(y)y_j,
\]
the only factor dependent on $h$ being $D_ju(x)$.
Recall $|\Jq(u(x)-u(y))|=|u(x)-u(y)|^{q-1}$ and, since $\Psi_\lam(z):=\psi_\lam(z)z$
is Lipschitz with $\operatorname{supp}\Psi_\lam\subseteq B^{\lam}:=\{|z|\le2/\lam\}$,
\[
  |w_j(x,y)|\le L_\lam\min\{|x-y|,C_\lam\}\,\mathbf 1_{E^{\lam}}(x,y)\,,
  \quad E^{\lam}:=\{x\in B^{\lam}\ \text{or}\ y\in B^{\lam}\},
\]
$L_\lam=\|\Psi_\lam\|_{\mathrm{Lip}}$, $C_\lam=2\,\mathrm{diam}\,B^{\lam}$.

\medskip

Define
\begin{equation}\label{eq:B1limqlep}
  B_1^{\lim}:=\sum_{j=1}^N\iint_{\R^{2N}}\Jq(u(x)-u(y))\,w_j(x,y)\,u_{x_j}(x)\,\dmu .
\end{equation}

\textbf{Claim.}
The integrand in \eqref{eq:B1limqlep} is
absolutely $\dmu$-integrable, $B_1^{\lim}$ is finite, and
$\displaystyle\lim_{h\to0}B_1=B_1^{\lim}$. We split the proof of the claim in three cases as follows:

\paragraph{Case 1: $q< p$.}

Set $v_h:=D_ju-u_{x_j}$, so $\|v_h\|_{L^p(\RN)}\to0$ by \eqref{eq:Djconv}. Then
\begin{equation}\label{eq:B1diffqlep}
  B_1(h)-B_1^{\lim}
  =\sum_{j=1}^N\iint_{\R^{2N}}\Jq(u(x)-u(y))\,w_j(x,y)\,v_h(x)\,\dmu .
\end{equation}
Applying H\"older's inequality on $(\R^{2N},\dmu)$ with exponents $\big(\tfrac{q}{q-1},q\big)$ gives
\begin{equation}\label{eq:B1diffholderqlep}
  \big|B_1(h)-B_1^{\lim}\big|
  \le[u]_{s,q}^{\,q-1}\sum_{j=1}^N\big(K_2(h)\big)^{1/q},
  \qquad
  K_2(h):=\iint_{\R^{2N}}|w_j|^{q}\,|v_h(x)|^{q}\,\dmu .
\end{equation}

By Tonelli, $K_2(h)=\int_{\RN}|v_h(x)|^q\,G(x)\,dx$ with
$G(x):=\int_{\RN}\frac{|w_j(x,y)|^{q}}{|x-y|^{N+sq}}\,dy$. Further, using $|w_j|\le L_\lam|x-y|$ on
$\{|x-y|<1\}$ and $|w_j|\le L_\lam C_\lam$ on $\{|x-y|\ge1\}$,
\[
  G(x)\le L_\lam^{q}\!\int_{|x-y|<1}\!\!|x-y|^{q-N-sq}dy
   + L_\lam^{q} C_\lam^{q}\!\int_{|x-y|\ge1}\!\!|x-y|^{-N-sq}dy
   =\frac{L_\lam^{q}|\mathbb{S}^{N-1}|}{q(1-s)}+\frac{L_\lam^{q} C_\lam^{q}|\mathbb{S}^{N-1}|}{sq}
   < \infty .
\]
Thus $\|G\|_\infty<\infty$. For $|x|>2R_\lam$, where $R_\lam:=2/\lam,\;w_j\ne0$ implies $y\in B^{\lam}$,
then $|x-y|\ge|x|/2$ and $G(x)\le C_\lam'\,|x|^{-(N+sq)}$.

Split $$K_2(h)=\int_{|x|\le2R_\lam}|v_h|^q G\,dx+\int_{|x|>2R_\lam}|v_h|^q G\,dx.$$
By H\"older with exponent $\tfrac pq>1$,
\[
  \int_{|x|\le2R_\lam}|v_h|^q\,G\,dx
  \le\|G\|_\infty\,|B_{2R_\lam}|^{1-\frac qp}\,\|v_h\|_{p}^{q}
  \xrightarrow[h\to0]{}0 .
\]
And
\[
  \int_{|x|>2R_\lam}|v_h|^q\,G\,dx
  \le C_\lam'\,\|v_h\|_{p}^{q}
   \Big(\int_{|x|>2R_\lam}|x|^{-(N+sq)(p/q)'}dx\Big)^{1/(p/q)'}
  \xrightarrow[h\to0]{}0,
\]
finite since $(N+sq)(p/q)'-N=\tfrac{q(N+ps)}{p-q}>0$. Hence $K_2(h)\to0$.

Replacing $v_h$ by $u_{x_j}\in L^p(\RN)$ in the above bounds gives
$|B_1^{\lim}|\le[u]_{s,q}^{q-1}\sum_j\big(\int_{\RN}|u_{x_j}|^q G\big)^{1/q}
<\infty$, so the integrand in \eqref{eq:B1limqlep} is absolutely
$\dmu$-integrable. With \eqref{eq:B1diffholderqlep}, $\lim_{h\to0}B_1=B_1^{\lim}$.

\paragraph{Case 2: $sq<p\leq q$.}
In this case apply H\"older on $(\R^{2N},\dmu)$ with exponents $p',p$ in \eqref{eq:B1diffqlep}, we obtain
\begin{equation}\label{eq:B1diffholder}
  \big|B_1(h)-B_1^{\lim}\big|
  \le\sum_{j=1}^N(K_1')^{1/p'}\,(K_2'(h))^{1/p},\qquad
  \begin{aligned}
  K_1'&:=\iint|u(x)-u(y)|^{(q-1)p'}\,\dmu,\\
  K_2'(h)&:=\iint|w_j|^{p}\,|v_h(x)|^{p}\,\dmu .
  \end{aligned}
\end{equation}

Since $q\geq p$, $(q-1)p'-q=\tfrac{q-p}{p-1}\ge0$, so by $u\in L^\infty(\rnn)$,
\[
  |u(x)-u(y)|^{(q-1)p'}\le(2M)^{(q-1)p'-q}|u(x)-u(y)|^{q},\; \text{ where } \;M:=\|u\|_\infty,
\]
so $K_1'\le(2M)^{(q-1)p'-q}\,[u]_{s,q}^q<\infty$.
$K_2'(h)=\int_{\RN}|v_h(x)|^p\,G_p(x)\,dx$ with
$G_p(x):=\int_{\RN}\frac{|w_j(x,y)|^{p}}{|x-y|^{N+sq}}\,dy$. Then similarly as in Case 1, we have 
$\|G_p\|_\infty<\infty$ as $p>sq$. Thus
\[K_2'(h)\le\|G_p\|_\infty\,\|v_h\|_{L^p(\RN)}^p\xrightarrow[h\to0]{}0 .
\]
Replacing $v_h$ by $u_{x_j}\in L^p(\RN)$ gives $|B_1^{\lim}|\le\sum_j(K_1')^{1/p'}
\big(\|G_p\|_\infty\|u_{x_j}\|_{p}^p\big)^{1/p}<\infty$, so the integrand in
\eqref{eq:B1limqlep} is absolutely $\dmu$-integrable. With
\eqref{eq:B1diffholder}, $\lim_{h\to0}B_1=B_1^{\lim}$.

\paragraph{Case 3: $p\leq sq \,(<q)$.}
\[
  H_j(x):=\int_{\RN}\Jq(u(x)-u(y))\,w_j(x,y)\,\frac{dy}{|x-y|^{N+sq}},
\]
for which $B_1(h)=\sum_j\int_{\RN}H_j\,D_ju\,dx$ and
$B_1^{\lim}=\sum_j\int_{\RN}H_j\,u_{x_j}\,dx$, showing $H_j\in L^{p'}(\RN)$; the
limit then follows from
\[
  \Big|B_1(h)-B_1^{\lim}\Big|
  =\Big|\sum_{j=1}^N\int_{\RN}H_j\,(D_ju-u_{x_j})\,dx\Big|
  \le\sum_{j=1}^N\|H_j\|_{{p'}}\,\|v_h\|_{p}
  \xrightarrow[h\to0]{}0,
\]
by \eqref{eq:Djconv}. In this case, by the assumption of Theorem~\ref{thm:main Pohozaev}, $u\in C^{0,l}_{\mathrm{loc}}(\RN)$ for some $l\in(l_*,1]$. Hence
\begin{equation}\label{eq:alphaassump}
  u\in C^{0,l}(B_\lam^{+2})
  \quad\text{with}\quad
  l>l_*=\frac{sq-p}{q-p},
\end{equation}
\noindent
where
$B_\lam^{+r}:=\{x:\operatorname{dist}(x,B^\lam)\le r\}$ for $r>0$. 

\medskip
Split $H_j=H_j^{\mathrm{near}}+H_j^{\mathrm{far}}$ at $|x-y|=1$. If $|x-y|<1$ and
$x\notin B_\lam^{+1}$ then $y\notin B^{\lam}$, so $w_j=0$. Thus
$\operatorname{supp}H_j^{\mathrm{near}}\subseteq B_\lam^{+1}$. By
$|w_j|\le L_\lam|x-y|$,
\begin{equation}\label{eq:Hnearcase3}
  |H_j^{\mathrm{near}}(x)|
  \le L_\lam\int_{|x-y|<1}\frac{|u(x)-u(y)|^{q-1}}{|x-y|^{N+sq-1}}\,dy .
\end{equation}
Set $t:=\dfrac{q-p}{p}\in(0,q-1)$. For $x\in B_\lam^{+1}$, $|x-y|<1$ both points lie in
$B_\lam^{+2}$, so \eqref{eq:alphaassump} gives
$|u(x)-u(y)|^{t}\le[u]_{C^{0,l}(B_\lam^{+2})}^{t}|x-y|^{l t}$.
Inserting this in \eqref{eq:Hnearcase3} and applying H\"older in $y$ with exponents $\big(\tfrac{q}{q-1-t},\tfrac{q}{1+t}\big)= (\frac{p}{p-1}, p)$, we obtain
\begin{equation}\label{eq:Hnearinterp}
\begin{aligned}
|H_j^{\mathrm{near}}(x)| &\le L_\lam\int_{|x-y|<1}\frac{|u(x)-u(y)|^{q-1-t}}{|x-y|^{(N+sq)\frac{q-1-t}{q}}}\frac{|u(x) - u(y)|^t}{|x-y|^{-1+\frac{(N+sq)(1+t)}{q}}}\,dy \\
& \leq C\, L_\lam \left( \int_{|x-y|<1}\frac{|u(x)-u(y)|^q}{|x-y|^{N+sq}}\,dy \right)^{\frac{q-1-t}{q}} \Big(\int_{|x-y|<1}|x-y|^{(l t+1)\frac{q}{1+t}-(N+sq)}\,dy\Big)^{\frac{1+t}{q}}
\end{aligned}
\end{equation}
where the second H\"older factor is a finite constant because in polar coordinates, the radial exponent exceeds $-1$ iff $\frac{l t+1}{1+t}>s$, i.e.\ $t(s-l)<1-s$, which for $t=\frac{q-p}{p}$ is exactly $l>l_*$, and this holds by \eqref{eq:alphaassump}. So, by \eqref{eq:Hnearinterp} and $\operatorname{supp}H_j^{\mathrm{near}}
\subseteq B_\lam^{+1}$,
\[
  \int_{\RN}|H_j^{\mathrm{near}}|^{p'}
  =\int_{B_\lam^{+1}}|H_j^{\mathrm{near}}|^{p'}
  \le C\,[u]_{s,q}^q<\infty .
\]

\medskip
Now we prove $H_j^{\mathrm{far}}\in L^{p'}$.
For $|x-y|\ge1$,
$|H_j^{\mathrm{far}}(x)|\le(2M)^{q-1}L_\lam C_\lam\int_{|z|\ge1}|z|^{-N-sq}dz=:K_{\mathrm{far}}$
for all $x$. And for $|x|>2R_\lam$, $w_j\ne0$ forces $y\in B^{\lam}$, whence
$|x-y|\ge|x|/2$ and $|H_j^{\mathrm{far}}(x)|\le C_\lam'|x|^{-(N+sq)}$. Therefore
\[
  \int_{\RN}|H_j^{\mathrm{far}}|^{p'}
  \le K_{\mathrm{far}}^{p'}|B_{2R_\lam}|
   +(C_\lam')^{p'}\!\int_{|x|>2R_\lam}\!\!|x|^{-(N+sq)p'}dx<\infty, \text{ since $(N+sq)p'>N$ }.
\]

\medskip
Thus $\lim_{h\to0}B_1=B_1^{\lim}$ in all the cases. 

\medskip
We now rewrite $B_1^{\lim}$ without derivatives of $u$. Fix $y\in\RN$ and $0<\varepsilon<R$, and set
\[
G_y(x):=\tfrac1q|u(x)-u(y)|^q,\qquad K_j(x):=w_j(x,y)\,|x-y|^{-N-sq}.
\]
Since $u\in W^{1,p}_{\mathrm{loc}}(\RN)\cap C(\RN)\cap L^\infty(\RN)$, we have $G_y\in W^{1,1}_{\mathrm{loc}}(\RN)\cap C(\RN)$ with $\partial_{x_j}G_y=\Jq(u(x)-u(y))\,u_{x_j}(x)$, while $K_j$ is Lipschitz on $\overline{B_R(y)}\setminus B_\varepsilon(y)$. The divergence theorem on $B_R(y)\setminus B_\varepsilon(y)$ gives
\begin{equation*}
\begin{aligned}
  \int_{B_R(y)\setminus B_\varepsilon(y)}(\partial_{x_j}G_y)\,K_j\,dx
  &=-\int_{B_R(y)\setminus B_\varepsilon(y)}G_y\,(\partial_{x_j}K_j)\,dx
   +\int_{\partial B_\varepsilon(y)}G_y\,K_j\,\nu_j\,d\sigma(x) \\
 & \quad   +\int_{\partial B_R(y)}G_y\,K_j\,\nu_j\,d\sigma(x),
   \end{aligned}
\end{equation*}
where $\nu=(\nu_1,\dots,\nu_N)$ is the outward unit normal to $B_R(y)\setminus B_\varepsilon(y)$. On $\partial B_R(y)$ we have $|w_j(x,y)|\le|\psi_\lam(x)x_j|
+|\psi_\lam(y)y_j|\le2\|\Psi_\lam\|_\infty$, then $|K_j(x)|\le
2\|\Psi_\lam\|_\infty R^{-N-sq}$, and together with $|G_y(x)|\le\tfrac1q(2\|u\|_\infty)^q$ this yields
\[
  \Big|\int_{\partial B_R(y)}G_y\,K_j\,\nu_j\,d\sigma\Big|
  \le\frac{2\|\Psi_\lam\|_\infty}{q}(2\|u\|_\infty)^q\,|\mathbb{S}^{N-1}|\,R^{-1-sq}
  \xrightarrow[R\to\infty]{}0 .
\]
Moreover,
\[
  \sum_{j=1}^N\partial_{x_j}K_j
  =\frac{\dive\Psi_\lam(x)}{|x-y|^{N+sq}}
   -(N+sq)\frac{(x-y)\cdot(\Psi_\lam(x)-\Psi_\lam(y))}{|x-y|^{N+sq+2}} .
\]
For a.e. $y$ all the integrals over $\RN\setminus B_\varepsilon(y)$ below are absolutely convergent (by the Claim and since $|\dive\Psi_\lam|\le C$ and $|(x-y)\cdot(\Psi_\lam(x)-\Psi_\lam(y))|\le L_\lam|x-y|^2$). Letting $R\to\infty$, summing over $j$ and integrating in $y$, we obtain for every $\varepsilon>0$
\begin{equation}\label{eq:B1mu}
\begin{aligned}
&\sum_{j=1}^N\iint_{\{|x-y|>\varepsilon\}}\Jq(u(x)-u(y))\,w_j(x,y)\,u_{x_j}(x)\,\dmu\\
&\quad=-\frac1q\iint_{\{|x-y|>\varepsilon\}}|u(x)-u(y)|^q\,\dive\Psi_\lam(x)\,\dmu\\
&\qquad+\frac{N+sq}{q}\iint_{\{|x-y|>\varepsilon\}}|u(x)-u(y)|^q\,
       \frac{(x-y)\cdot(\Psi_\lam(x)-\Psi_\lam(y))}{|x-y|^{2}}\,\dmu
+I_\varepsilon,
\end{aligned}
\end{equation}
where
\[
I_\varepsilon:=\frac1q\int_{\RN}\int_{\partial B_\varepsilon(y)}|u(x)-u(y)|^q\,\frac{(\Psi_\lam(x)-\Psi_\lam(y))\cdot\nu(x)}{|x-y|^{N+sq}}\,d\sigma(x)\,dy .
\]
As $\varepsilon\to0^+$, the left-hand side of \eqref{eq:B1mu} tends to $B_1^{\lim}$ by the absolute integrability of its integrand, and the first two terms on the right-hand side tend to the corresponding integrals over $\R^{2N}$ by dominated convergence, since $|u(x)-u(y)|^q\in L^1(\R^{2N},\dmu)$.

It remains to show that $I_{\varepsilon_k}\to0$ along some sequence $\varepsilon_k\to0^+$. On $\partial B_\varepsilon(y)$, we have $|x-y|=\varepsilon$ and $|\Psi_\lam(x)-\Psi_\lam(y)|\le L_\lam\varepsilon$, so
\[
  |I_\varepsilon|
  \le\frac{L_\lam}{q}\,\varepsilon^{\,1-N-sq}\,\Upsilon(\varepsilon),
  \qquad
  \Upsilon(\varepsilon):=\int_{\RN}\int_{\partial B_\varepsilon(y)}|u(x)-u(y)|^q\,d\sigma(x)\,dy .
\]
In polar coordinates centered at $y$,
\[
  [u]_{s,q}^q=\iint_{\R^{2N}}\frac{|u(x)-u(y)|^q}{|x-y|^{N+sq}}\,dx\,dy
  =\int_0^\infty \Upsilon(\varepsilon)\,\varepsilon^{-N-sq}\,d\varepsilon,
\]
so the nonnegative function $\upsilon(\varepsilon):=\Upsilon(\varepsilon)\,\varepsilon^{-N-sq}$ belongs to $L^1(0,1)$, and $|I_\varepsilon|\le\frac{L_\lam}{q}\,\varepsilon\,\upsilon(\varepsilon)$. Since $\upsilon\in L^1(0,1)$, we have $\liminf_{\varepsilon\to0^+}\varepsilon\,\upsilon(\varepsilon)=0$, otherwise $\varepsilon\upsilon(\varepsilon)\ge c>0$ for small $\varepsilon$ would force
$\upsilon(\varepsilon)\ge c/\varepsilon$, which is not integrable near $0$. Hence there is a sequence $\varepsilon_k\to0^+$ with $\varepsilon_k\upsilon(\varepsilon_k)\to0$, so $I_{\varepsilon_k}\to0$. Passing to the limit in \eqref{eq:B1mu} along $\varepsilon_k$, we conclude
\begin{equation}\label{eq:B1final}
\begin{aligned}
\lim_{h\to0}B_1=B_1^{\lim}
&=-\frac1q\iint_{\R^{2N}}|u(x)-u(y)|^q\,\dive\Psi_\lam(x)\,\dmu\\
&\quad+\frac{N+sq}{q}\iint_{\R^{2N}}|u(x)-u(y)|^q\,
       \frac{(x-y)\cdot(\Psi_\lam(x)-\Psi_\lam(y))}{|x-y|^{N+sq+2}}\,dx\,dy .
\end{aligned}
\end{equation}
Note that no H\"older continuity of $u$ has been used in this step.

\paragraph{Limit of $B_{2,1}$.}
Apply the nonlocal discrete integration by parts \eqref{eq:dibp-nl} of
Lemma~\ref{lem:null Identity}\,(b), and then the change of variables $(x,y)\mapsto(x-he_j,y-he_j)$, which leaves $\dmu$ invariant
\begin{align*}
B_{2,1}&=\frac1q\sum_{j=1}^N\iint_{\R^{2N}}D_j^{\mathrm{diag}}\big(|u(x)-u(y)|^q\big)\,
     \psi_\lam(x)x_j\dmu \\
 &=-\frac1q\sum_{j=1}^N\iint_{\R^{2N}}
     |u(x+he_j)-u(y+he_j)|^q\,D_j(\psi_\lam(x)x_j)\dmu\\
 &=-\frac1q\sum_{j=1}^N\iint_{\R^{2N}}
     |u(x)-u(y)|^q\,\frac{\psi_{\lambda}(x) x_j - \psi_{\lambda}(x-he_j)(x_j - h)}{h}\dmu.
\end{align*}
The last difference quotient is bounded by the Lipschitz constant of $x\mapsto\psi_\lam(x)x_j$, uniformly in $h$, and converges to $\partial_{x_j}(\psi_\lam(x)x_j)$ uniformly on $\RN$. Since $|u(x)-u(y)|^q\in L^1(\R^{2N},\dmu)$, the dominated convergence theorem gives
\begin{align*}
  \lim_{h\to0}B_{2,1}
  &=-\frac1q\sum_{j=1}^N\iint_{\R^{2N}}
      |u(x)-u(y)|^q\,\partial_{x_j}(\psi_\lam(x)x_j)\dmu \notag\\
  &=-\frac Nq\iint_{\R^{2N}}|u(x)-u(y)|^q\psi_\lam(x)\dmu
    -\frac1q\iint_{\R^{2N}}|u(x)-u(y)|^q\,\nabla\psi_\lam(x)\cdot x\dmu.
\end{align*}

Combining the estimates from the local and nonlocal parts (see \eqref{eq:B1final}), we obtain
\begin{align}
  \lim_{h\to0}(A+B)
  &=\int_{\RN}|\nabla u|^{p-2}(\nabla u\cdot\nabla\psi_\lam)(x\cdot\nabla u)\,dx
   +\int_{\RN}\psi_\lam|\nabla u|^p\,dx \notag\\
  &\quad-\frac Np\int_{\RN}|\nabla u|^p\psi_\lam\,dx
   -\frac1p\int_{\RN}|\nabla u|^p\nabla\psi_\lam\cdot x\,dx\notag \\
  &\quad-\frac Nq\iint_{\R^{2N}}|u(x)-u(y)|^q\psi_\lam(x)\dmu
    -\frac1q\iint_{\R^{2N}}|u(x)-u(y)|^q\,\nabla\psi_\lam(x)\cdot x\dmu \notag\\
  &\quad-\frac1q\iint_{\R^{2N}}
       |u(x)-u(y)|^q \,\dive\Psi_\lam(x)\dmu \notag\\
   &\quad +\frac{(N+sq)}{q}\iint_{\R^{2N}}
       |u(x)-u(y)|^q\frac{(x-y)\cdot(\Psi_\lam(x)-\Psi_\lam(y))}
       {|x-y|^{N+sq+2}}\,dx\,dy. \label{eq:ABlimit}
\end{align}
\subsubsection{The right-hand side L}
\[
L=L(h)=\int_{\RN}f(u)\,\psi_\lam\sum_{j=1}^N x_j\,D_ju\,dx
   =\sum_{j=1}^N\int_{\RN}E_j(x)\,D_ju(x)\,dx,
   \quad E_j:=f(u)\,\psi_\lam\,x_j.
\]
We pass to the limit by H\"older
exactly as for $A_1,A_2$. Indeed $E_j\in L^{p'}(\RN)$ follows from $f(u)\in L^{p'}(\RN)$, see \eqref{eq:integrability}, and the boundedness of $\psi_\lam x_j$. Thus
\begin{align}
  &\lim_{h\to0}L
  =\sum_{j=1}^N\int_{\RN}f(u)\psi_\lam x_j\,u_{x_j}\,dx \notag \\
&=\int_{\RN}\psi_\lam\,x\cdot\nabla\big(F(u)\big)\,dx = -\Big(N\int_{\RN}\psi_\lam F(u)\,dx
+\int_{\RN}F(u)\,x\cdot\nabla\psi_\lam\,dx\Big). \label{eq:Llim} \end{align}
In the second equality, we used the chain rule \eqref{eq:chain}. In the last equality, we integrated by parts, which is allowed since $F(u)\in W^{1,1}_{\mathrm{loc}}(\RN)$ and the vector field $\psi_\lam x$ is of class $C^1$ with compact support, and $\dive(\psi_\lam x)=N\psi_\lam+x\cdot\nabla\psi_\lam$.

\subsubsection{Proof of Theorem~\ref{thm:main Pohozaev}}
Insert \eqref{eq:ABlimit}, \eqref{eq:Llim} into the $h\to0$ limit
of \eqref{eq:master}, i.e.\ $\lim_h A+\lim_h B=\lim_h L$,
\begin{align*}
  &-\Big(N\!\int_{\RN}\!\psi_\lam F(u)\,dx
        +\int_{\RN}\!F(u)\,x\cdot\nabla\psi_\lam\,dx\Big) \notag\\
  &=\int_{\RN}|\nabla u|^{p-2}(\nabla u\cdot\nabla\psi_\lam)(x\cdot\nabla u)\,dx
        +\frac{p-N}{p}\int_{\RN}|\nabla u|^p\psi_\lam\,dx
        -\frac1p\int_{\RN}|\nabla u|^p\,\nabla\psi_\lam\cdot x\,dx \notag\\
  &\qquad-\frac Nq\iint_{\R^{2N}}|u(x)-u(y)|^q\psi_\lam(x)\dmu
        -\frac1q\iint_{\R^{2N}}|u(x)-u(y)|^q\,\nabla\psi_\lam(x)\cdot x\dmu \notag\\
  &\qquad-\frac1q\iint_{\R^{2N}}
       \frac{|u(x)-u(y)|^q}{|x-y|^{N+sq}}\,\dive\Psi_\lam(x)\,dx\,dy \notag\\
  &\qquad+\frac{(N+sq)}{q}\iint_{\R^{2N}}
       |u(x)-u(y)|^q\frac{(x-y)\cdot(\Psi_\lam(x)-\Psi_\lam(y))}
       {|x-y|^{N+sq+2}}\,dx\,dy.
\end{align*}
Now we let $\lam\to0$. As $\lam\to0$ we have, for all $x,y\in\RN$ with $x\neq y$,
\[
\psi_\lam(x)\to1,\qquad \nabla\psi_\lam(x)\cdot x\to0,\qquad \dive\Psi_\lam(x)=N\psi_\lam(x)+\nabla\psi_\lam(x)\cdot x\to N,
\]
\[
\frac{(x-y)\cdot(\Psi_\lam(x)-\Psi_\lam(y))}{|x-y|^2}\to1 .
\]
Moreover, by \eqref{eq:cutoff}, all these quantities are bounded independently of $\lam$. In particular, $|\Psi_\lam(x)-\Psi_\lam(y)|\le C|x-y|$ with $C$ independent of $\lam$. Since $|\nabla u|^p\in L^1(\RN)$, $F(u)\in L^1(\RN)$, $|u(x)-u(y)|^q\in L^1(\R^{2N},\dmu)$ and
$\big||\nabla u|^{p-2}(\nabla u\cdot\nabla\psi_\lam)(x\cdot\nabla u)\big|\le|x||\nabla\psi_\lam||\nabla u|^p\le C|\nabla u|^p$, the dominated convergence theorem allows us to pass to the limit in every term. The first term on the right-hand side and the terms containing $\nabla\psi_\lam\cdot x$ tend to $0$, and we obtain
\[
  -N\int_{\RN}F(u)\,dx
  =\Big(1-\frac Np\Big)[u]_{1,p}^p+\Big(-\frac Nq-\frac Nq+\frac{N+sq}{q}\Big)[u]_{s,q}^q
  =-\frac{N-p}{p}[u]_{1,p}^p-\frac{N-sq}{q}[u]_{s,q}^q .
\]
The proof of Theorem~\ref{thm:main Pohozaev} is thus complete.\hfill$\square$

\section{Nonexistence}\label{sec:nonex}
In this section, we combine Theorem~\ref{thm:main Pohozaev} with the Nehari identity. No relation between $\ps$ and $q^*_s$ is assumed. Recall that $m_-=\min\{\ps,q^*_s\}$ and $m_+=\max\{\ps,q^*_s\}$.
Let $u$ be as in Theorem~\ref{thm:main Pohozaev}, and assume in addition that $f(u)u\in L^1(\RN)$. Then $\varphi=u$ is admissible in \eqref{eq:weak}, and we obtain the Nehari identity
\begin{equation}\label{eq:nehari}
  [u]_{1,p}^p+[u]_{s,q}^q=\int_{\RN}f(u)\,u\,dx .
\end{equation}
Since $\frac{N-p}{p}=\frac{N}{\ps}$ and $\frac{N-sq}{q}=\frac{N}{q^*_s}$, \eqref{eq:poho} gives
\begin{equation}\label{eq:pohoN}
  \frac{1}{\ps}[u]_{1,p}^p+\frac{1}{q^*_s}[u]_{s,q}^q=\int_{\RN}F(u)\,dx .
\end{equation}
Subtracting $\tau$ times \eqref{eq:nehari} from \eqref{eq:pohoN}, we get, for every $\tau\in\R$,
\begin{equation}\label{eq:keyrel}
  \Big(\frac{1}{\ps}-\tau\Big)[u]_{1,p}^p
  +\Big(\frac{1}{q^*_s}-\tau\Big)[u]_{s,q}^q
  =\int_{\RN}\Big(F(u)-\tau\, f(u)\,u\Big)\,dx .
\end{equation}
We use \eqref{eq:keyrel} with $\tau=1/m_+$ and with $\tau=1/m_-$. For $\tau=1/m_+$ both coefficients on the left-hand side are nonnegative, and for $\tau=1/m_-$ both are nonpositive. If $\ps\neq q^*_s$, then in each case exactly one of the two coefficients is zero, and if $\ps=q^*_s$, then both coefficients are zero.

\begin{proof}[Proof of Theorem \ref{thm:nonex}]
Case (i). Take $\tau=1/m_+$ in \eqref{eq:keyrel}. By (i), $F(u)-\frac{1}{m_+}f(u)u\le0$ in $\RN$, so the right-hand side of \eqref{eq:keyrel} is $\le0$, while the left-hand side is $\ge0$. Hence both sides are zero.

If $\ps<q^*_s$, then $m_+=q^*_s$ and the coefficient $\frac1{\ps}-\frac1{q^*_s}$ of $[u]_{1,p}^p$ is positive, so $[u]_{1,p}=0$. Similarly, if $\ps>q^*_s$, then $[u]_{s,q}=0$. By the Sobolev inequalities $\|u\|_{L^{\ps}(\RN)}\le C[u]_{1,p}$ and $\|u\|_{L^{q^*_s}(\RN)}\le C[u]_{s,q}$, in both cases $u=0$ a.e.\ in $\RN$, hence $u\equiv0$ by continuity.

If $\ps=q^*_s$, the inequality in (i) is strict for $t\neq0$. Then the function $F(u)-\frac{1}{m_+}f(u)u$ is $\le0$ in $\RN$ and $<0$ on the set $\{u\neq0\}$. Since its integral is zero, the set $\{u\neq0\}$ has measure zero, and again $u\equiv0$.

Case (ii) is treated in the same way, with $\tau=1/m_-$ and all inequalities reversed. 
\end{proof}

\section{Existence in the intermediate range}\label{sec:existence}

In Section~\ref{sec:nonex}, we saw that, for the power nonlinearity in \eqref{eq:Pr}, the Pohozaev and Nehari identities leave open only the exponents $r$ with $m_-<r<m_+$. The aim of this section is to prove Theorem \ref{thm:existenceNonexistence}, Theorem \ref{thm:groundstate} and Theorem \ref{thm:same}. Throughout the section, we assume \eqref{assump: standing} and $\ps\neq q^*_s$, so that $m_-<m_+$.

\subsection{Functional setting}\label{sec:functional}

Recall that $\mathcal S$ and $\mathcal S_f$ denote the best constants in the Sobolev inequalities
\[
\mathcal S\,\|u\|_{\ps}^p\le[u]_{1,p}^p,\qquad \mathcal S_f\,\|u\|_{q^*_s}^q\le[u]_{s,q}^q,\qquad u\in\X(\RN).
\]
We recall that, since $1<p<N$ and $sq<N$, we have $D^{1,p}_0(\RN)=\{u\in L^{\ps}(\RN):\nabla u\in L^p(\RN)\}$ and $D^{s,q}_0(\RN)=\{u\in L^{q^*_s}(\RN):[u]_{s,q}<\infty\}$. The space $\X(\mathbb{R}^N)=D^{1,p}_0(\RN)\cap D^{s,q}_0(\RN)$ is reflexive. Indeed, the map
\begin{equation}
T u := \left(\nabla u, \frac{u(x)-u(y)}{|x-y|^{(N+sq)/q}}\right)\label{eq:isometryMap}
\end{equation}
is an isometry from $\X(\mathbb{R}^N)$ into the product space $L^p(\mathbb{R}^N;\mathbb{R}^N)\times L^q(\mathbb{R}^{2N})$ endowed with the norm
\[
\|(F, G)\|_{L^p \times L^q} := \|F\|_{L^p(\mathbb{R}^N;\mathbb{R}^N)} + \|G\|_{L^q(\mathbb{R}^{2N})},
\]
and its image is closed.

Let $$\X^{\mathrm{rad}}(\mathbb{R}^N) := \{ u \in \X(\mathbb{R}^N) : u \text{ is radial}\}.$$ It is a closed linear subspace of $\X(\mathbb{R}^N)$, and consequently, it is weakly closed.

For $r>1$, we recall from \eqref{eq:JG-intro} the energy functional $\mathcal J$ and the Pohozaev functional $\mathcal G$ associated with \eqref{eq:Pr}, and the Pohozaev set $\mathcal M$. We also set
\[
\Phi(u):=\frac1p[u]_{1,p}^p+\frac1q[u]_{s,q}^q ,
\]
so that $\mathcal J(u)=\Phi(u)-\frac1r\|u\|_r^r$, and $\mathcal S_r=\inf\{\Phi(u):\ u\in\X(\RN),\ \|u\|_r=1\}$ whenever $m_-\le r\le m_+$.
By Lemma~\ref{lem:interp} below, $\mathcal J$ and $\mathcal G$ are well defined and of class $C^1$ on $\X(\RN)$ whenever $m_-\le r\le m_+$. Since $F(t)=\frac1r|t|^r$ for $f(t)=|t|^{r-2}t$, a nontrivial weak solution of \eqref{eq:Pr} satisfies the Pohozaev identity \eqref{eq:poho} if and only if it belongs to $\mathcal M$. In particular, every nontrivial weak solution of \eqref{eq:Pr} to which Theorem~\ref{thm:main Pohozaev} applies lies in $\mathcal M$.

For $u\in\X(\RN)$ and $t>0$, as in the Introduction, $u_t(x)=u(x/t)$. Then
\begin{equation}\label{eq:dilation-scaling}
[u_t]_{1,p}^p=t^{N-p}[u]_{1,p}^p,\qquad [u_t]_{s,q}^q=t^{N-sq}[u]_{s,q}^q,\qquad \|u_t\|_r^r=t^N\|u\|_r^r .
\end{equation}
The map $t\mapsto u_t$ is continuous from $(0,\infty)$ into $\X(\RN)$, because dilations act continuously on $L^p(\RN;\RN)$ and on $L^q(\R^{2N})$. Moreover, $\|u_t\|_{\X}\to0$ as $t\to0^+$ by \eqref{eq:dilation-scaling}.

The first step is to find the exponents $r$ for which $\X(\RN)\subset L^r(\RN)$.

\begin{lemma}\label{lem:interp}
Let $m_-\le r\le m_+$, and let $\theta\in[0,1]$ be defined by $\frac1r=\frac{\theta}{\ps}+\frac{1-\theta}{q^*_s}$. Then, for every $u\in\X(\RN)$,
\begin{equation}\label{eq:interp}
\|u\|_r^r\le K\,\big([u]_{1,p}^p\big)^{\frac{\theta r}{p}}\big([u]_{s,q}^q\big)^{\frac{(1-\theta)r}{q}},
\qquad K:=\mathcal S^{-\frac{\theta r}{p}}\,\mathcal S_f^{-\frac{(1-\theta)r}{q}}>0 .
\end{equation}
Moreover,
\begin{equation}\label{eq:Lambda}
\Lambda:=\frac{\theta r}{p}+\frac{(1-\theta)r}{q}=1+\frac rN\big(\theta+(1-\theta)s\big)>1 ,
\end{equation}
and consequently, for every $u\in\X(\RN)$,
\begin{equation}\label{eq:interp-sum}
\|u\|_r^r\le K\big([u]_{1,p}^p+[u]_{s,q}^q\big)^{\Lambda}\qquad\text{and}\qquad \|u\|_r\le K^{\frac1r}\|u\|_{\X}.
\end{equation}
\end{lemma}

\begin{proof}
Since $u\in L^{\ps}(\RN)\cap L^{q^*_s}(\RN)$, H\"older's inequality gives $\|u\|_r\le\|u\|_{\ps}^{\theta}\|u\|_{q^*_s}^{1-\theta}$, and \eqref{eq:interp} follows from the two Sobolev inequalities. Using $\frac1p=\frac1{\ps}+\frac1N$ and $\frac1q=\frac1{q^*_s}+\frac sN$ gives \eqref{eq:Lambda}. The first inequality in \eqref{eq:interp-sum} follows from \eqref{eq:interp}, because each of $[u]_{1,p}^p$ and $[u]_{s,q}^q$ is bounded by their sum. The second one again follows from \eqref{eq:interp} and the inequality $a^{\theta}b^{1-\theta}\le a+b$ for $a,b\ge0$. \qedhere
\end{proof}

\subsection{The compact embedding}
We shall use two classical tools in this section. The first one is the radial lemma for $D^{1,p}_0(\mathbb{R}^N)$ (see \cite{berestycki_lions_1983, su_wang_willem_2007}): since $1 < p < N$, every radially symmetric $u \in D^{1,p}_0(\mathbb{R}^N)$ satisfies
\begin{equation}\label{eq:radial-lemma}
|u(x)| \le C(N,p)\,[u]_{1,p}\,|x|^{-\frac{N-p}{p}} \qquad \text{for a.e. } x \in \mathbb{R}^N.
\end{equation}

The second tool is the following compactness lemma of Strauss type (see \cite{strauss_1977} and the appendix of \cite{berestycki_lions_1983}). 

\begin{lemma}\label{lem:LDCT_Dr}
Let $P, Q : \mathbb{R} \to \mathbb{R}$ be continuous functions with $\frac{P(t)}{Q(t)} \to 0$ as $|t| \to \infty$ and as $t \to 0$. Let $(w_n)$ be a sequence of measurable functions such that 
\[
\sup_n \int_{\mathbb{R}^N} |Q(w_n)|\,dx < \infty, \quad P(w_n) \to g \text{ a.e. in } \mathbb{R}^N,
\]
and $w_n(x) \to 0$ as $|x| \to \infty$, uniformly in $n$. Then $P(w_n) \to g$ in $L^1(\mathbb{R}^N)$.
\end{lemma}
\begin{proposition}\label{prop:embedding}
Let $r>1$.
\begin{enumerate}[label=\textup{(\alph*)}]
\item The embedding $\X(\RN)\hookrightarrow L^r(\RN)$ holds if and only if $m_-\le r\le m_+$.
\item The embedding $\X^{\mathrm{rad}}(\RN)\hookrightarrow L^r(\RN)$ is compact if and only if $m_-<r<m_+$.
\end{enumerate}
\end{proposition}

\begin{proof}
If $m_-\le r\le m_+$, the embedding in (a) is \eqref{eq:interp-sum}. Conversely, if $\X(\RN)\subset L^r(\RN)$, then the inclusion map has a closed graph, since convergence in either space implies a.e.\ convergence along a subsequence. By the closed graph theorem, it is continuous. It is therefore enough to show that the inequality $\|u\|_r\le C\|u\|_{\X}$ fails when $r\notin[m_-,m_+]$, and that the embedding in (b) is not compact when $r=m_\pm$.

For this, we fix a nonnegative radially symmetric $\varphi\in C^\infty_c(\RN)\setminus\{0\}$, and for $a>0$ and $\lam>0$ we set $\varphi^{a,\lam}(x):=\lam^a\varphi(\lam x)$. Then
\[
[\varphi^{a,\lam}]_{1,p}^p=\lam^{ap+p-N}[\varphi]_{1,p}^p,\qquad
[\varphi^{a,\lam}]_{s,q}^q=\lam^{aq+sq-N}[\varphi]_{s,q}^q,\qquad
\|\varphi^{a,\lam}\|_r^r=\lam^{ar-N}\|\varphi\|_r^r .
\]
We first take $a=N/m_+$. Then $ap+p-N=pN\big(\frac1{m_+}-\frac1{\ps}\big)\le0$ and $aq+sq-N=qN\big(\frac1{m_+}-\frac1{q^*_s}\big)\le0$, so $\|\varphi^{a,\lam}\|_{\X}\le\|\varphi\|_{\X}$ for $\lam\ge1$, while $\|\varphi^{a,\lam}\|_r^r=\lam^{N(r/m_+-1)}\|\varphi\|_r^r$. If $r>m_+$, this tends to $+\infty$ as $\lam\to\infty$, and the embedding in (a) fails. \\
If $r = m_+$, the $L^{m_+}$-norm of $\varphi^{a,\lambda}$ is invariant under $\lambda$. Since the family is bounded in $\X^{\mathrm{rad}}(\mathbb{R}^N)$ and $\varphi^{a,\lambda}(x) \to 0$ a.e.\ as $\lambda \to \infty$, we have $\varphi^{a,\lambda} \rightharpoonup 0$ weakly in $\X(\mathbb{R}^N)$. However, because $\|\varphi^{a,\lambda}\|_{m_+} = \|\varphi\|_{m_+} \neq 0$, the sequence cannot converge strongly to $0$ in $L^{m_+}(\mathbb{R}^N)$, thereby proving that the embedding in (b) is not compact.

\smallskip
\noindent
The exponent $m_-$ is treated in the same way, with $a=N/m_-$ and $\lam\to0^+$.

\smallskip
\noindent
To prove the compactness in (b), let $m_-<r<m_+$, and $(u_n)$ be bounded in $\X^{\mathrm{rad}}(\RN)$. Since $\X(\RN)$ is reflexive, up to a subsequence $u_n\rightharpoonup u$ in $\X(\RN)$ for some $u\in\X^{\mathrm{rad}}(\RN)$. For every $R>0$ the sequence $(u_n)$ is bounded in $W^{1,p}(B_R)$, because $\|u_n\|_{L^p(B_R)}\le|B_R|^{\frac1p-\frac1{\ps}}\|u_n\|_{\ps}$. By the Rellich--Kondrachov theorem and a diagonal argument, a further subsequence converges in $L^p_{\mathrm{loc}}(\RN)$ and a.e.\ in $\RN$, and its limit is $u$ because $u_n\rightharpoonup u$ in $L^{\ps}(\RN)$. We apply Lemma \ref{lem:LDCT_Dr} to $w_n:=u_n-u$, with $P(t)=|t|^r$ and $Q(t)=|t|^{\ps}+|t|^{q^*_s}$. Since $m_-<r<m_+$, we have $P(t)/Q(t)\le|t|^{r-m_-}\to0$ as $t\to0$ and $P(t)/Q(t)\le|t|^{r-m_+}\to0$ as $|t|\to\infty$. Moreover, $\sup_n\int_{\RN}Q(w_n)\,dx<\infty$ by the Sobolev inequalities, and $P(w_n)\to0$ a.e. Finally, $w_n$ is radially symmetric and $\sup_n[w_n]_{1,p}<\infty$, so \eqref{eq:radial-lemma} gives $\sup_n|w_n(x)|\le C|x|^{-\frac{N-p}{p}}\to0$ as $|x|\to\infty$. Therefore $\|u_n-u\|_r^r=\int_{\RN}P(w_n)\,dx\to0$.
\end{proof}

\subsection{Existence via quotient minimization: Proof of Theorem \ref{thm:existenceNonexistence}}
Minimizing the energy under the constraint $\|u\|_r=1$ produces a weak solution of \eqref{eq:Pr} with a Lagrange multiplier $\mu > 0$ in front of $|u|^{r-2}u$. Unlike single-operator problems, $\mu$ cannot be removed by a simple scaling because \eqref{eq:Pr} lacks scaling invariance in the subcritical-supercritical regime $m_- < r < m_+$.

Indeed, a rescaling $u \mapsto \kappa\,u(\lam\,\cdot)$ for $\kappa, \lam > 0$ multiplies $[u]_{1,p}^p$, $[u]_{s,q}^q$, and $\|u\|_r^r$ by $\kappa^p\lam^{p-N}$, $\kappa^q\lam^{sq-N}$, and $\kappa^r\lam^{-N}$, respectively. Leaving the functional $\mathcal J$ invariant requires $\frac{N-p}{p}=\frac{N-sq}{q}=\frac{N}{r}$, i.e., $\ps = q^*_s = r$, whereas leaving the equation \eqref{eq:Pr} alone invariant requires $D(r)=0$, with $D(r)$ as in \eqref{eq:Dr}, which occurs only when $p > sq$ and $r = {\tilde r} := \frac{pq(1-s)}{p-sq}$ with $\lam = \kappa^{\beta}$ and $\beta = \frac{q-p}{p-sq}$.

As shown in Lemma~\ref{lem:Dnonzero} below, $D(r) \neq 0$ throughout $[m_-, m_+]$ whenever $\ps \neq q^*_s$. Consequently, while no one-parameter scaling leaves \eqref{eq:Pr} invariant, we can still absorb the Lagrange multiplier $\mu$ by performing a simultaneous, decoupled change of amplitude and spatial scale performed in Lemma~\ref{lem:normalisation}.

\begin{lemma}\label{lem:Dnonzero}
Assume $\ps\neq q^*_s$ and $m_-\le r\le m_+$. Then $D(r)\neq0$.
\end{lemma}

\begin{proof}
If $p\le sq$, then $D(r)\le-pq(1-s)<0$. Now let $p>sq$, and recall $\tilde{r}=\frac{pq(1-s)}{p-sq}$, so that $D(r)=(p-sq)(r-\tilde{r})$. Setting $\delta:=N(q-p)-pq(1-s)$, a direct computation gives
\[
q^*_s-\ps=\frac{N\,\delta}{(N-p)(N-sq)},\qquad
\tilde{r}-q^*_s=\frac{sq\,\delta}{(p-sq)(N-sq)} .
\]
Since $\ps\neq q^*_s$, we have $\delta\neq0$, and both differences share the sign of $\delta$. Hence either $\ps<q^*_s<\tilde{r}$ or $\tilde{r}<q^*_s<\ps$. In both cases $\tilde{r}\notin[m_-,m_+]$, and therefore $D(r)\neq0$.
\end{proof}

\begin{lemma}\label{lem:normalisation}
Assume $\ps\neq q^*_s$ and $m_-\le r\le m_+$, let $\mu>0$, and let $u\in\X(\RN)$ satisfy \eqref{eq:Pmu}. Then the function $v(x)=\kappa\,u(\lam x)$ defined in \eqref{eq:normalisation} is a weak solution of \eqref{eq:Pr} in the sense of Definition~\ref{def:weak}. 
\end{lemma}
\begin{proof}
The exponents in \eqref{eq:normalisation} are the unique solution $(a,b)$ of the linear system $(p-q)a+(p-sq)b=0$, $(p-r)a+pb=-1$, whose determinant is $D(r)$ ($\neq 0$ by Lemma~\ref{lem:Dnonzero}). Hence $\kappa=\mu^a$ and $\lam=\mu^b$\; for $a= \frac{p-sq}{D(r)},\; b= \frac{q-p}{D(r)}$ satisfy
\begin{equation}\label{eq:normalisation-system}
\kappa^{p-1}\lam^{p}=\kappa^{q-1}\lam^{sq}=\kappa^{r-1}\mu^{-1}=:\Theta .
\end{equation}
Let $\varphi\in\X(\RN)$ and set $\psi(y):=\varphi(y/\lam)$, which again belongs to $\X(\RN)$. All the integrals below are finite by Lemma~\ref{lem:interp}. The change of variables $y=\lam x$ gives
\[
\langle v,\varphi\rangle_p=\kappa^{p-1}\lam^{p-N}\langle u,\psi\rangle_p,\qquad
\langle v,\varphi\rangle_{s,q}=\kappa^{q-1}\lam^{sq-N}\langle u,\psi\rangle_{s,q},
\]
and
\[
\int_{\RN}|v|^{r-2}v\,\varphi\,dx=\kappa^{r-1}\lam^{-N}\int_{\RN}|u|^{r-2}u\,\psi\,dx .
\]
By \eqref{eq:normalisation-system}, and \eqref{eq:Pmu} tested with $\psi$,
\[
\begin{aligned}
\langle v,\varphi\rangle_p+\langle v,\varphi\rangle_{s,q}
&=\lam^{-N}\Theta\big(\langle u,\psi\rangle_p+\langle u,\psi\rangle_{s,q}\big)
=\lam^{-N}\Theta \mu\int_{\RN}|u|^{r-2}u \psi dx
=\int_{\RN}|v|^{r-2}v \varphi dx .
\end{aligned}
\]
This is the weak formulation of \eqref{eq:Pr}.
\end{proof}

We conclude this section by proving Theorem \ref{thm:existenceNonexistence} via standard arguments.

\begin{proof}[Proof of Theorem \ref{thm:existenceNonexistence}]
For (a) and (b) we use, as in the proof of Proposition~\ref{prop:embedding}, the functions $\varphi^{a,\lam}(x)=\lam^a\varphi(\lam x)$ with $\varphi\in C^\infty_c(\RN)\setminus\{0\}$. Recall that
\[
[\varphi^{a,\lam}]_{1,p}^p=\lam^{ap+p-N}[\varphi]_{1,p}^p,\qquad
[\varphi^{a,\lam}]_{s,q}^q=\lam^{aq+sq-N}[\varphi]_{s,q}^q,\qquad
\|\varphi^{a,\lam}\|_r^r=\lam^{ar-N}\|\varphi\|_r^r .
\]

(a)Since $C^\infty_c(\RN)\subset\X(\RN)\cap L^r(\RN)$, the set over which the infimum is taken is
nonempty for every $r>1$. Normalize $\|\varphi\|_r=1$ and take $a=N/r$, so that $\|\varphi^{a,\lam}\|_r=1$ for every $\lam>0$. The two energy exponents are
\[
ap+p-N=pN\Big(\frac1r-\frac1{\ps}\Big),\qquad aq+sq-N=qN\Big(\frac1r-\frac1{q^*_s}\Big).
\]
If $r>m_+$, both are negative, and the energy of $\varphi^{a,\lam}$ tends to $0$ as $\lam\to\infty$. If $r<m_-$, both are positive, and the energy tends to $0$ as $\lam\to0^+$. In both cases $\mathcal S_r=0$.

\medskip
(b) Since $\X(\RN)\subset L^{p^*}(\RN)\cap L^{q^*_s}(\RN)$, for $r\in\{p^*,q^*_s\}$ the infimum is
taken over all $u\in\X(\RN)$ with $\|u\|_r=1$. Let $r=\ps$. For every $u\in\X(\RN)$ with $\|u\|_{\ps}=1$, the Sobolev inequality gives
\[
\frac1p[u]_{1,p}^p+\frac1q[u]_{s,q}^q\ge\frac1p\,\mathcal S+\frac1q[u]_{s,q}^q>\frac1p\,\mathcal S ,
\]
since $[u]_{s,q}>0$ for $u\neq0$. Hence $\mathcal S_r\ge\frac1p\mathcal S$, and the value $\frac1p\mathcal S$ cannot be attained.\\ For the opposite inequality, let $\varepsilon>0$. Since $C^\infty_c(\RN)$ is dense in $D^{1,p}_0(\RN)$, there is $\varphi\in C^\infty_c(\RN)$ with $\|\varphi\|_{\ps}=1$ and $[\varphi]_{1,p}^p\le\mathcal S+\varepsilon$. Take $a=N/\ps$. Then $\|\varphi^{a,\lam}\|_{\ps}=1$ and $[\varphi^{a,\lam}]_{1,p}=[\varphi]_{1,p}$ for every $\lam>0$, while
$$[\varphi^{a,\lam}]_{s,q}^q=\lam^{qN\left(\frac1{\ps}-\frac1{q^*_s}\right)}[\varphi]_{s,q}^q .
$$
This tends to $0$ as $\lam\to0^+$ if $\ps<q^*_s$, and as $\lam\to\infty$ if $\ps>q^*_s$. Hence $\mathcal S_r\le\frac1p(\mathcal S+\varepsilon)$, and since $\varepsilon$ is arbitrary, $\mathcal S_r=\frac1p\mathcal S$. The case $r=q^*_s$ follows analogously.

\medskip
For part (c), we proceed by direct minimization. If $\|u\|_r=1$, then \eqref{eq:interp-sum} gives $[u]_{1,p}^p+[u]_{s,q}^q\ge K^{-1/\Lambda}$, hence $\mathcal S_r\ge\frac{1}{\max\{p,q\}}K^{-1/\Lambda}>0$.

Let $(u_n)$ be a minimizing sequence. We have $[|u_n|]_{1,p}=[u_n]_{1,p}$ and $[|u_n|]_{s,q}\le[u_n]_{s,q}$. Let $u_n^\star$ be the Schwarz symmetrization of $|u_n|$, which is well defined since $|u_n|\in L^{\ps}(\RN)$. By the P\'olya--Szeg\H{o} inequality and its nonlocal counterpart (see \cite{almgren_lieb_1989, ourPreviousPaper_isocriticalregimemixedlocalnonlocal}), $[u_n^\star]_{1,p}\le[u_n]_{1,p}$ and $[u_n^\star]_{s,q}\le[u_n]_{s,q}$. Moreover, $u_n^\star\in L^{\ps}(\RN)\cap L^{q^*_s}(\RN)$ and $\|u_n^\star\|_r=\|u_n\|_r=1$. Hence $u_n^\star\in\X^{\mathrm{rad}}(\RN)$, and $(u_n^\star)$ is again a minimizing sequence. We may therefore assume that every $u_n$ is nonnegative, radially symmetric and radially nonincreasing. The sequence is bounded in $\X(\RN)$, since $[u_n]_{1,p}^p+[u_n]_{s,q}^q\le\max\{p,q\}(\mathcal S_r+1)$ for $n$ large.

\smallskip
By Proposition~\ref{prop:embedding}(b), up to a subsequence, $u_n\rightharpoonup u$ in $\X(\RN)$, $u_n\to u$ in $L^r(\RN)$ and a.e.\ in $\RN$. In particular $\|u\|_r=1$, and $u$ is nonnegative and radially symmetric. Moreover, $u$ agrees a.e.\ with a radially nonincreasing function. The functionals $u\mapsto[u]_{1,p}^p$ and $u\mapsto[u]_{s,q}^q$ are convex and continuous on $\X(\RN)$, hence weakly lower semicontinuous. Therefore, $\mathcal S_r$ is attained by $u$.

\smallskip
The functional $\Phi$ of Subsection~\ref{sec:functional} is of class $C^1$ on $\X(\RN)$, with $\langle\Phi'(u),\varphi\rangle=\langle u,\varphi\rangle_p+\langle u,\varphi\rangle_{s,q}$. By \eqref{eq:interp-sum}, the same holds for $\Psi(u):=\|u\|_r^r$, with $\langle\Psi'(u),\varphi\rangle=r\int_{\RN}|u|^{r-2}u\,\varphi\,dx$. Since $\langle\Psi'(u),u\rangle=r\neq0$ on the constraint, the Lagrange multiplier rule gives $\nu\in\R$ with $\Phi'(u)=\nu\Psi'(u)$. This is \eqref{eq:Pmu} with $\mu:=\nu r$. Testing with $\varphi=u$ gives $\mu=[u]_{1,p}^p+[u]_{s,q}^q>0$. Lemma~\ref{lem:normalisation} now shows that $v$ is a weak solution of \eqref{eq:Pr}. Since $\kappa,\lam>0$, $v$ has the same qualitative properties as $u$.
\end{proof}
\subsection{The fibering map}
Since $r$ may be smaller than $\max\{p,q\}$, the function $t \mapsto \mathcal{J}(tu)$ can tend to $+\infty$ as $t \to \infty$ rather than attaining a global maximum, which makes the Nehari manifold unsuitable as a constraint. To obtain a unique maximum along one-parameter paths and capture ground states, we instead consider $\mathcal{J}$ along dilations $t \mapsto u_t$ as defined in \eqref{eq:dilation-scaling}.
\begin{lemma}\label{lem:fibering}
Let $m_-<r<m_+$, and let $u\in\X(\RN)$ with $\|u\|_r>0$. Then the function $h_u(t):=\mathcal J(u_t)$, $t>0$, has exactly one critical point $t_u>0$, which is its strict global maximum point, and $u_t\in\mathcal M$ if and only if $t=t_u$. Moreover,
\begin{enumerate}[label=\textup{(\roman*)}]
\item if $u\in\mathcal M$, then $t_u=1$ and
\begin{equation}\label{eq:J-on-M}
\mathcal J(u)=\max_{t>0}\mathcal J(u_t)=\frac1N\Big([u]_{1,p}^p+s[u]_{s,q}^q\Big);
\end{equation}
\item $c_{\mathcal M}=\inf_{\mathcal M}\mathcal J>0$.
\end{enumerate}
\end{lemma}

\begin{proof}
By \eqref{eq:dilation-scaling},
\[
h_u(t)=\frac{t^{N-p}}{p}[u]_{1,p}^p+\frac{t^{N-sq}}{q}[u]_{s,q}^q-\frac{t^N}{r}\|u\|_r^r\,,
\qquad
t\,h_u'(t)=\mathcal G(u_t).
\]
Since $u\neq0$, both $[u]_{1,p}$ and $[u]_{s,q}$ are positive, so the function
\[
t^{-N}\,\mathcal G(u_t)=\frac{N-p}{p}\,t^{-p}[u]_{1,p}^p+\frac{N-sq}{q}\,t^{-sq}[u]_{s,q}^q-\frac Nr\|u\|_r^r
\]
is strictly decreasing on $(0,\infty)$, and it tends to $+\infty$ as $t\to0^+$ and to $-\frac Nr\|u\|_r^r<0$ as $t\to\infty$. Hence it has exactly one zero $t_u$, the function $h_u$ is increasing on $(0,t_u)$ and decreasing on $(t_u,\infty)$, and $u_t\in\mathcal M$ if and only if $t=t_u$.

If $u\in\mathcal M$, then $\frac Nr\|u\|_r^r=\frac{N-p}{p}[u]_{1,p}^p+\frac{N-sq}{q}[u]_{s,q}^q>0$, and $t_u=1$ because $\mathcal G(u_1)=\mathcal G(u)=0$. The identity in \eqref{eq:J-on-M} follows from
\[
\mathcal J(u)=\mathcal J(u)-\frac1N\,\mathcal G(u)=\Big(\frac1p-\frac{N-p}{Np}\Big)[u]_{1,p}^p+\Big(\frac1q-\frac{N-sq}{Nq}\Big)[u]_{s,q}^q .
\]
To prove (ii), let $u\in\mathcal M$ and set $\gamma_0:=\min\{\frac{N-p}{p},\frac{N-sq}{q}\}>0$. By \eqref{eq:interp-sum},
\[
\gamma_0\big([u]_{1,p}^p+[u]_{s,q}^q\big)\le\frac Nr\|u\|_r^r\le\frac{NK}{r}\big([u]_{1,p}^p+[u]_{s,q}^q\big)^{\Lambda},
\]
and since $\Lambda>1$, we obtain $[u]_{1,p}^p+[u]_{s,q}^q\ge\rho_0:=\big(\frac{\gamma_0 r}{NK}\big)^{\frac1{\Lambda-1}}>0$. Then \eqref{eq:J-on-M} gives $c_{\mathcal M}\ge\frac sN\rho_0>0$.
\end{proof}

\subsection{Ground states: Proof of Theorem \ref{thm:groundstate}}

In this subsection, we prove Theorem \ref{thm:groundstate}. To obtain a nonnegative solution, we work with the truncated functionals
\[
\mathcal J_+(u):=\frac1p[u]_{1,p}^p+\frac1q[u]_{s,q}^q-\frac1r\|u^+\|_r^r,\qquad
\mathcal G_+(u):=\frac{N-p}{p}[u]_{1,p}^p+\frac{N-sq}{q}[u]_{s,q}^q-\frac Nr\|u^+\|_r^r,
\]
which are of class $C^1$ on $\X(\RN)$, and coincide with $\mathcal J$ and $\mathcal G$ on nonnegative functions. Let
\begin{equation}\label{eq:cmp}
\Gamma:=\big\{\gamma\in C([0,1],\X^{\mathrm{rad}}(\RN)):\ \gamma(0)=0,\ \mathcal J_+(\gamma(1))<0\big\},\qquad
c_{\mathrm{mp}}:=\inf_{\gamma\in\Gamma}\max_{\tau\in[0,1]}\mathcal J_+(\gamma(\tau)) .
\end{equation}

\begin{lemma}\label{lem:jeanjean}
Let $m_-<r<m_+$. Then $\Gamma\neq\emptyset$ and $c_{\mathrm{mp}}>0$. Moreover, there exists a sequence $(u_n)\subset\X^{\mathrm{rad}}(\RN)$ such that
\begin{equation}\label{eq:PS-Pohozaev}
\mathcal J_+(u_n)\to c_{\mathrm{mp}},\qquad \|\mathcal J_+'(u_n)\|_{(\X^{\mathrm{rad}}(\RN))^*}\to0,\qquad \mathcal G_+(u_n)\to0 .
\end{equation}
\end{lemma}

\begin{proof}
Let $\varphi\in C^\infty_c(\RN)\setminus\{0\}$ be nonnegative and radially symmetric. By the explicit form of $h_\varphi$ in the proof of Lemma~\ref{lem:fibering}, $\mathcal J_+(\varphi_t)=h_\varphi(t)\to-\infty$ as $t\to\infty$, because $N>\max\{N-p,N-sq\}$. Since $\varphi_t\to0$ in $\X(\RN)$ as $t\to0^+$, the path $\tau\mapsto\varphi_{\tau t_0}$, extended by $0$ at $\tau=0$, belongs to $\Gamma$ for $t_0$ large. So $\Gamma \neq \emptyset$.\\
Next, for every $u\in\X(\RN)$, since $\|u^+\|_r\le\|u\|_r$, \eqref{eq:interp-sum} gives
\[
\mathcal J_+(u)\ge g\big([u]_{1,p}^p+[u]_{s,q}^q\big),\qquad g(\tau):=\frac{\tau}{\max\{p,q\}}-\frac Kr\,\tau^{\Lambda}.
\]
Because $\Lambda>1$, there is $\rho_1>0$ with $g>0$ on $(0,\rho_1]$. If $\gamma\in\Gamma$, then $\mathcal J_+(\gamma(1))<0$ forces $[\gamma(1)]_{1,p}^p+[\gamma(1)]_{s,q}^q>\rho_1$. By continuity, $\gamma$ meets the set $\{[u]_{1,p}^p+[u]_{s,q}^q=\rho_1\}$, on which $\mathcal J_+\ge g(\rho_1)$. Hence $c_{\mathrm{mp}}\ge g(\rho_1)>0$.

\smallskip
We now consider the Banach space $E:=\R\times\X^{\mathrm{rad}}(\RN)$, with the norm $|\sigma|+\|v\|_{\X}$, and the functional $\widetilde{\mathcal J}(\sigma,v):=\mathcal J_+(v_{e^{\sigma}})$, where $v_{e^\sigma}(x)=v(e^{-\sigma}x)$. By \eqref{eq:dilation-scaling},
\[
\widetilde{\mathcal J}(\sigma,v)=\frac{e^{(N-p)\sigma}}{p}[v]_{1,p}^p+\frac{e^{(N-sq)\sigma}}{q}[v]_{s,q}^q-\frac{e^{N\sigma}}{r}\|v^+\|_r^r ,
\]
so $\widetilde{\mathcal J}\in C^1(E,\R)$, and
\begin{equation}\label{eq:augmented-derivatives}
\partial_\sigma\widetilde{\mathcal J}(\sigma,v)=\mathcal G_+(v_{e^\sigma}),\qquad
\langle\partial_v\widetilde{\mathcal J}(\sigma,v),\psi\rangle=\langle\mathcal J_+'(v_{e^\sigma}),\psi_{e^\sigma}\rangle\quad\text{for }\psi\in\X^{\mathrm{rad}}(\RN).
\end{equation}
The second identity follows from the chain rule, since $v\mapsto v_{e^\sigma}$ is linear and continuous on $\X^{\mathrm{rad}}(\RN)$.

Let $\widetilde\Gamma$ be the set of $\widetilde\gamma\in C([0,1],E)$ with $\widetilde\gamma(0)=(0,0)$ and $\widetilde{\mathcal J}(\widetilde\gamma(1))<0$, and let $\widetilde c$ be the corresponding minimax value. If $\gamma\in\Gamma$, then $(0,\gamma)\in\widetilde\Gamma$, so $\widetilde c\le c_{\mathrm{mp}}$. Conversely, if $\widetilde\gamma=(\sigma,v)\in\widetilde\Gamma$, then $\gamma(\tau):=v(\tau)_{e^{\sigma(\tau)}}$ belongs to $\Gamma$, because the map $(\sigma,v)\mapsto v_{e^\sigma}$ is continuous from $E$ into $\X^{\mathrm{rad}}(\RN)$ (by \eqref{eq:dilation-scaling}, $\|v_{e^\sigma}-w_{e^{\sigma'}}\|_{\X}\le e^{N|\sigma|}\|v-w\|_{\X}+\|w_{e^\sigma}-w_{e^{\sigma'}}\|_{\X}$). Since $\mathcal J_+(\gamma(\tau))=\widetilde{\mathcal J}(\widetilde\gamma(\tau))$, we get $\widetilde c\ge c_{\mathrm{mp}}$. Thus $\widetilde c=c_{\mathrm{mp}}>0$, whereas $\max\{\widetilde{\mathcal J}(0,0),\widetilde{\mathcal J}(\widetilde\gamma(1))\}=0$ for every $\widetilde\gamma\in\widetilde\Gamma$.

For $n\in\mathbb N$ large, choose $\gamma_n\in\Gamma$ with $\max_{[0,1]}\mathcal J_+\circ\gamma_n\le c_{\mathrm{mp}}+\frac1{n^2}$. The general minimax principle \cite[Theorem~2.8]{willem_minimax_1996}, applied to $\widetilde{\mathcal J}$, $\widetilde\Gamma$, the path $(0,\gamma_n)$, $\varepsilon=\frac1{n^2}$ and $\delta=\frac1n$, gives $(\sigma_n,v_n)\in E$ with
\[
\big|\widetilde{\mathcal J}(\sigma_n,v_n)-c_{\mathrm{mp}}\big|\le\frac{2}{n^2},\qquad
\mathrm{dist}\big((\sigma_n,v_n),\{0\}\times\gamma_n([0,1])\big)\le\frac2n,\qquad
\|\widetilde{\mathcal J}'(\sigma_n,v_n)\|_{E^*}\le\frac8n .
\]
In particular $|\sigma_n|\le\frac2n$. Set $u_n:=(v_n)_{e^{\sigma_n}}\in\X^{\mathrm{rad}}(\RN)$. Then $\mathcal J_+(u_n)=\widetilde{\mathcal J}(\sigma_n,v_n)\to c_{\mathrm{mp}}$, and $\mathcal G_+(u_n)=\partial_\sigma\widetilde{\mathcal J}(\sigma_n,v_n)\to0$ by \eqref{eq:augmented-derivatives}. Let $\varphi\in\X^{\mathrm{rad}}(\RN)$ and $\psi:=\varphi(e^{\sigma_n}\cdot)$, so that $\psi_{e^{\sigma_n}}=\varphi$. Since $[\psi]_{1,p}=e^{\sigma_n\frac{p-N}{p}}[\varphi]_{1,p}$ and $[\psi]_{s,q}=e^{\sigma_n\frac{sq-N}{q}}[\varphi]_{s,q}$, we have $\|\psi\|_{\X}\le e^{N|\sigma_n|}\|\varphi\|_{\X}$, and \eqref{eq:augmented-derivatives} gives
\[
|\langle\mathcal J_+'(u_n),\varphi\rangle|=|\langle\partial_v\widetilde{\mathcal J}(\sigma_n,v_n),\psi\rangle|\le\frac8n\,e^{\frac{2N}{n}}\,\|\varphi\|_{\X}.
\]
This proves \eqref{eq:PS-Pohozaev}.
\end{proof}

With the Palais--Smale sequence provided by Lemma~\ref{lem:jeanjean} in hand, we are now ready to establish the existence of a radial ground state.
\begin{proof}[Proof of Theorem \ref{thm:groundstate}]
Let $(u_n)$ be given by Lemma~\ref{lem:jeanjean}. As in the proof of Lemma~\ref{lem:fibering}, $\mathcal J_+-\frac1N\mathcal G_+=\frac1N\big([\,\cdot\,]_{1,p}^p+s[\,\cdot\,]_{s,q}^q\big)$, so
\[
\frac1N\Big([u_n]_{1,p}^p+s[u_n]_{s,q}^q\Big)=\mathcal J_+(u_n)-\frac1N\mathcal G_+(u_n)\to c_{\mathrm{mp}} ,
\]
and $(u_n)$ is bounded in $\X(\RN)$. By Proposition~\ref{prop:embedding}(b), up to a subsequence, $u_n\rightharpoonup w$ in $\X(\RN)$ and $u_n\to w$ in $L^r(\RN)$, for some $w\in\X^{\mathrm{rad}}(\RN)$.

We claim that $u_n\to w$ in $\X(\RN)$. Since $(u_n-w)$ is bounded, $\langle\mathcal J_+'(u_n),u_n-w\rangle\to0$ by \eqref{eq:PS-Pohozaev}, and $\langle\mathcal J_+'(w),u_n-w\rangle\to0$ by weak convergence. Moreover, by H\"older's inequality,
\[
\Big|\int_{\RN}\big((u_n^+)^{r-1}-(w^+)^{r-1}\big)(u_n-w)\,dx\Big|\le\big(\|u_n\|_r^{r-1}+\|w\|_r^{r-1}\big)\|u_n-w\|_r\to0 .
\]
Consequently,
\[
\big(\langle u_n,u_n-w\rangle_p-\langle w,u_n-w\rangle_p\big)+\big(\langle u_n,u_n-w\rangle_{s,q}-\langle w,u_n-w\rangle_{s,q}\big)\to0 .
\]
By H\"older's inequality in $\RN$,
\[
\langle u_n,u_n-w\rangle_p-\langle w,u_n-w\rangle_p\ge\big([u_n]_{1,p}^{p-1}-[w]_{1,p}^{p-1}\big)\big([u_n]_{1,p}-[w]_{1,p}\big)\ge0,
\]
and in the same way, using H\"older's inequality in $(\R^{2N},\dmu)$, the nonlocal term is bounded below by $\big([u_n]_{s,q}^{q-1}-[w]_{s,q}^{q-1}\big)\big([u_n]_{s,q}-[w]_{s,q}\big)\ge0$.

Hence both terms tend to $0$, which gives $[u_n]_{1,p}\to[w]_{1,p}$ and $[u_n]_{s,q}\to[w]_{s,q}$. Since also $Tu_n\rightharpoonup Tw$, the uniform convexity of $L^p(\RN;\RN)$ and $L^q(\R^{2N})$ gives $Tu_n\to Tw$, where $T$ is as defined in \eqref{eq:isometryMap}. This proves the claim.

By continuity, $\mathcal J_+(w)=c_{\mathrm{mp}}>0$, $\mathcal G_+(w)=0$, and $\mathcal J_+'(w)=0$ in $(\X^{\mathrm{rad}}(\RN))^*$. In particular $w\neq0$. The orthogonal group $O(N)$ acts isometrically on $\X(\RN)$ by $(g,u)\mapsto u\circ g^{-1}$, $\mathcal J_+$ is invariant under this action, and $\X^{\mathrm{rad}}(\RN)$ is the set of fixed points. By the principle of symmetric criticality \cite{palais_1979}, $\mathcal J_+'(w)=0$ in $\X(\RN)^*$.

We test this equation with $w^-:=\max\{-w,0\}\in\X(\RN)$. Then the standard estimates give
\[
0=\langle\mathcal J_+'(w),w^-\rangle\le-[w^-]_{1,p}^p-[w^-]_{s,q}^q ,
\]
so $w^-=0$, that is, $w\ge0$. Therefore $\mathcal J'(w)=\mathcal J_+'(w)=0$. Thus $w$ is a nontrivial, nonnegative, radially symmetric weak solution of \eqref{eq:Pr}. Moreover $\mathcal G(w)=\mathcal G_+(w)=0$, so $w\in\mathcal M$ and $c_{\mathcal M}\le\mathcal J(w)=c_{\mathrm{mp}}$.

It remains to prove that $c_{\mathrm{mp}}\le c_{\mathcal M}$. Let $z\in\mathcal M$, and let $\tilde z$ be the Schwarz symmetrization of $|z|$. As in the proof of Theorem~\ref{thm:existenceNonexistence}, $\tilde z\in\X^{\mathrm{rad}}(\RN)$, $[\tilde z]_{1,p}\le[z]_{1,p}$, $[\tilde z]_{s,q}\le[z]_{s,q}$ and $\|\tilde z^+\|_r=\|z\|_r>0$. Hence, by \eqref{eq:dilation-scaling} and \eqref{eq:J-on-M},
\[
\mathcal J_+(\tilde z_t)\le\mathcal J(z_t)\le\mathcal J(z)\qquad\text{for every }t>0.
\]
Since $\mathcal J_+(\tilde z_t)\to-\infty$ as $t\to\infty$ and $\tilde z_t\to0$ in $\X(\RN)$ as $t\to0^+$, the path $\tau\mapsto\tilde z_{\tau t_0}$, extended by $0$ at $\tau=0$, belongs to $\Gamma$ for $t_0$ large. Therefore $c_{\mathrm{mp}}\le\sup_{t>0}\mathcal J_+(\tilde z_t)\le\mathcal J(z)$, and taking the infimum over $z\in\mathcal M$ gives $c_{\mathrm{mp}}\le c_{\mathcal M}$. Hence $\mathcal J(w)=c_{\mathcal M}=c_{\mathrm{mp}}$. \qedhere
\end{proof}

\subsection{Energy equivalence of two minimization problems}

Theorems~\ref{thm:existenceNonexistence} and~\ref{thm:groundstate} find solutions by two
different minimizations. The first Theorem gives $v$ which comes from the quotient $\mathcal S_r$, and whereas the second theorem gives $w$ which comes from the energy on the Pohozaev set. We now show that $v$ also has the least
energy, $\mathcal J(v)=c_{\mathcal M}$. Behind both problems lies the best
constant in the interpolation inequality \eqref{eq:interp}.

Let $m_-<r<m_+$ and let $\theta\in(0,1)$ be as in Lemma~\ref{lem:interp}, so that $\frac1r=\frac{\theta}{\ps}+\frac{1-\theta}{q^*_s}$. This gives
\begin{equation}\label{eq:theta-id}
\frac{1}{\ps}-\frac1r=(1-\theta)\,k,\qquad
\frac1r-\frac{1}{q^*_s}=\theta\,k,\qquad
k:=\frac{1}{\ps}-\frac{1}{q^*_s}\neq0 .
\end{equation}
We set
\begin{equation}\label{eq:Qdef}
\mathcal Q(u):=\frac{[u]_{1,p}^{\theta}\,[u]_{s,q}^{1-\theta}}{\|u\|_r},
\qquad
\mathcal C_r=\inf_{u\in\X(\RN)\setminus\{0\}}\mathcal Q(u)>0 ,
\end{equation}
where $\mathcal C_r$ is the constant defined in the Introduction, that is, the best constant in \eqref{eq:interp}. By the 
choice of $\theta$,
$$\mathcal Q(\kappa\,u(\lam\,\cdot\,))=\mathcal Q(u)\qquad\text{for all }\kappa>0,\ \lam>0.$$
 Finally, recall from the Introduction that
\[
\frac{1}{p_\theta}:=\frac{\theta}{p}+\frac{1-\theta}{q},
\qquad
c_\theta:=\theta^{\theta/p}(1-\theta)^{(1-\theta)/q} .
\]
Since $p<\ps$ and $q<q^*_s$, we have
$p_\theta<r$.

\begin{lemma}\label{lem:rigid}
Let $\mu>0$, and let $z\in\X(\RN)\setminus\{0\}$ satisfy \eqref{eq:Pmu} and
the Pohozaev identity
$\frac{N-p}{p}[z]_{1,p}^p+\frac{N-sq}{q}[z]_{s,q}^q=\frac{N\mu}{r}\|z\|_r^r$.
Then
\[
[z]_{1,p}^p=\theta\mu\|z\|_r^r,\qquad
[z]_{s,q}^q=(1-\theta)\mu\|z\|_r^r,\qquad
\mathcal Q(z)=c_\theta\,\mu^{1/p_\theta}\,\|z\|_r^{\,r/p_\theta-1}.
\]
\end{lemma}

\begin{proof}
Write $A=[z]_{1,p}^p$, $B=[z]_{s,q}^q$, $C=\|z\|_r^r$. Testing \eqref{eq:Pmu} with $z$ gives $A+B=\mu C$, and the Pohozaev identity reads $\frac{A}{\ps}+\frac{B}{q^*_s}=\frac{\mu C}{r}$. Then using \eqref{eq:theta-id}, we get $(1-\theta)A=\theta B$. Hence $A=\theta\mu C$ and $B=(1-\theta)\mu C$, and the formula for then $\mathcal Q(z)$ follows from \eqref{eq:Qdef}.
\end{proof}
A weak solution of \eqref{eq:Pr} satisfies \eqref{eq:Pmu} with $\mu=1$, since $|z|^{r-1}\varphi\in L^1(\RN)$ for all $\varphi\in\X(\RN)$ by Lemma~\ref{lem:interp}. So Lemma~\ref{lem:rigid} applies, with $\mu=1$, to every solution of \eqref{eq:Pr} in $\mathcal M$, and for such $z$
\begin{equation}\label{eq:J-rigid}
\mathcal J(z)=\frac{[z]_{1,p}^p}{p}+\frac{[z]_{s,q}^q}{q}-\frac{\|z\|_r^r}{r}
=\Big(\frac1{p_\theta}-\frac1r\Big)\|z\|_r^r
=\frac{r-p_\theta}{r\,p_\theta}\Big(\frac{\mathcal Q(z)}{c_\theta}\Big)^{\frac{r\,p_\theta}{r-p_\theta}} .
\end{equation}

\begin{lemma}\label{lem:minQ}
Let $u$, $\mu$ and $v$ be as in Theorem~\ref{thm:existenceNonexistence}. Then $u$ satisfies the Pohozaev identity of Lemma~\ref{lem:rigid}, $v\in\mathcal M$, and $\mathcal Q(v)=\mathcal Q(u)=\mathcal C_r$.
\end{lemma}
\begin{proof}
Let $\Phi$ be as in Subsection~\ref{sec:functional} and, for $\eta>0$, let $z^{(\eta)}:=\eta^{N/r}z(\eta\,\cdot\,)$ as in the Introduction. Note that, $\|z^{(\eta)}\|_r=\|z\|_r$ and
\begin{equation}\label{eq:Phi-eta}
\Phi(z^{(\eta)})=\alpha_1\eta^{a}+\alpha_2\eta^{b},
\qquad
\alpha_1:=\tfrac1p[z]_{1,p}^p,\quad \alpha_2:=\tfrac1q[z]_{s,q}^q,
\end{equation}
where, by \eqref{eq:theta-id},
\[
a:=pN\Big(\frac1r-\frac1{\ps}\Big)=-pN(1-\theta)k,
\qquad
b:=qN\Big(\frac1r-\frac1{q^*_s}\Big)=qN\theta k \qquad \text{and} \qquad ab <0.
\]

\smallskip
We first derive the identities satisfied by $u$ and $v$. Since $\|u^{(\eta)}\|_r=1$ for all $\eta>0$, each $u^{(\eta)}$ is admissible for $\mathcal S_r$, so $\eta\mapsto\Phi(u^{(\eta)})$ attains its minimum at $\eta=1$. By \eqref{eq:Phi-eta},
\[
0=\frac{d}{d\eta}\Phi(u^{(\eta)})\Big|_{\eta=1}=a\alpha_1+b\alpha_2
=N\Big(\frac1r-\frac1{\ps}\Big)[u]_{1,p}^p+N\Big(\frac1r-\frac1{q^*_s}\Big)[u]_{s,q}^q .
\]
Testing \eqref{eq:Pmu} with $\varphi=u$ gives $[u]_{1,p}^p+[u]_{s,q}^q=\mu\|u\|_r^r$, and therefore
\[
\frac{1}{\ps}[u]_{1,p}^p+\frac{1}{q^*_s}[u]_{s,q}^q=\frac{\mu}{r}\|u\|_r^r ,
\]
which is the Pohozaev identity of Lemma~\ref{lem:rigid}. By
\eqref{eq:normalisation-system},
\[
[v]_{1,p}^p=\kappa\Theta\lam^{-N}[u]_{1,p}^p,\qquad
[v]_{s,q}^q=\kappa\Theta\lam^{-N}[u]_{s,q}^q,\qquad
\|v\|_r^r=\kappa\Theta\lam^{-N}\mu\|u\|_r^r ,
\]
so $\mathcal G(v)=0$. Since $v\neq0$, we get $v\in\mathcal M$. Further, let
$\vartheta:=\frac{b}{b-a}=\frac{q\theta}{q\theta+p(1-\theta)}=\frac{\theta\,p_\theta}{p}\in(0,1)$,
so that $1-\vartheta=\frac{a}{a-b}$. The substitution
$\eta=(\alpha_2/\alpha_1)^{1/(a-b)}\xi$ gives
\[
\alpha_1\eta^a+\alpha_2\eta^b=\alpha_1^{\vartheta}\alpha_2^{1-\vartheta}\,(\xi^a+\xi^b),
\]
and therefore
\begin{equation}\label{eq:inf-eta}
\inf_{\eta>0}\big(\alpha_1\eta^a+\alpha_2\eta^b\big)=m_0\,\alpha_1^{\vartheta}\alpha_2^{1-\vartheta},
\qquad m_0:=\inf_{\xi>0}(\xi^a+\xi^b)\in[1,2],
\end{equation}
where $m_0\ge1$ because $ab<0$. If $\|z\|_r=1$, then, since
$p\vartheta=\theta p_\theta$ and $q(1-\vartheta)=(1-\theta)p_\theta$,
\[
\alpha_1^{\vartheta}\alpha_2^{1-\vartheta}
=p^{-\vartheta}q^{\vartheta-1}\big([z]_{1,p}^{\theta}[z]_{s,q}^{1-\theta}\big)^{p_\theta}
=p^{-\vartheta}q^{\vartheta-1}\,\mathcal Q(z)^{p_\theta}.
\]
Combining this with \eqref{eq:Phi-eta} and \eqref{eq:inf-eta}, and using that
each $z^{(\eta)}$ is admissible for $\mathcal S_r$, we obtain
\[
\mathcal S_r\le\inf_{\eta>0}\Phi(z^{(\eta)})=c_*\,\mathcal Q(z)^{p_\theta}
\quad\text{for every } z \text{ with } \|z\|_r=1,
\qquad c_*:=m_0\,p^{-\vartheta}q^{\vartheta-1}>0,
\]
with equality for $z=u$, because $\inf_{\eta}\Phi(u^{(\eta)})=\Phi(u)=\mathcal S_r$.
Hence $\mathcal Q(u)\le\mathcal Q(z)$ whenever $\|z\|_r=1$. Since $\mathcal Q$ is
invariant under the scalings $z\mapsto\kappa z(\lam\,\cdot\,)$, this gives
$\mathcal Q(u)=\mathcal C_r$, and also $\mathcal Q(v)=\mathcal Q(u)$.
\end{proof}

\begin{proof}[Proof of Theorem \ref{thm:same}]
Let $w$ be the ground state of Theorem~\ref{thm:groundstate}. Both $v$ and $w$ are solutions of \eqref{eq:Pr} lying in $\mathcal M$ (for $v$ by Lemma~\ref{lem:minQ}), so \eqref{eq:J-rigid} applies to both. Since $\mathcal Q(v)=\mathcal C_r\le\mathcal Q(w)$, we get $$\mathcal J(v)\le\mathcal J(w)=c_{\mathcal M}$$ because by \eqref{eq:J-rigid}, for the solutions lying in $\mathcal M$, the energy is an increasing function of $\mathcal Q$. And $\mathcal J(v)\ge c_{\mathcal M}$ because $v\in\mathcal M$. Hence $\mathcal J(v)=c_{\mathcal M}$, which also gives the second expression in \eqref{eq:cM-Sr}. For the first expression, we apply Lemma~\ref{lem:rigid} to $u$, for which $\|u\|_r=1$, and to $v$ with $\mu=1$. Since $\mathcal Q(v)=\mathcal Q(u)$ by Lemma \ref{lem:minQ}, we obtain
\[
\mathcal S_r=\Phi(u)=\frac{\theta\mu}{p}+\frac{(1-\theta)\mu}{q}=\frac{\mu}{p_\theta}
\qquad\text{and}\qquad
c_\theta\|v\|_r^{\,r/p_\theta-1}=\mathcal Q(v)=\mathcal Q(u)=c_\theta\,\mu^{1/p_\theta}.
\]
Hence $\|v\|_r^{\,r-p_\theta}=\mu=p_\theta\mathcal S_r$, and \eqref{eq:J-rigid} yields
\[
c_{\mathcal M}=\mathcal J(v)=\frac{r-p_\theta}{r\,p_\theta}\,\|v\|_r^{\,r}
=\frac{r-p_\theta}{r\,p_\theta}\,\big(p_\theta\mathcal S_r\big)^{\frac{r}{r-p_\theta}} . \qedhere\]
\end{proof}
\section*{Acknowledgements}
S.~Malhotra is supported by the Prime Minister's Research Fellows Fellowship (PMRF ID-1403266),
Government of India. Part of this work was carried out while D.~Gupta was a postdoctoral
researcher at KAIST, supported by
the Jang Young Sil Fellowship Program.
\section*{Declarations}
\textbf{Conflict of interest.} The authors declare that they have no conflict of interest.\\
\textbf{Data availability.} Data sharing is not applicable to this article, as no datasets were generated or analyzed during the current study.

\bibliographystyle{abbrv}
\bibliography{references}

\end{document}